\documentclass[11pt]{article}
\usepackage{amsmath, amssymb, amsfonts}
\usepackage{epsfig}

\def\be{\begin{equation}}
\def\ee{\end{equation}}
\def\bea{\begin{eqnarray}}
\def\eea{\end{eqnarray}}
\def\bes{\begin{eqnarray*}}
\def\ees{\end{eqnarray*}}

\def\nn{\nonumber}
\def\<{\langle}
\def\>{\rangle}
\def\lb{\label}
\def\bs{\setminus}

\def\R{{\mathbb{R}}}
\def\C{{\bf C}}
\def\Z{{\mathbb{Z}}}
\def\N{{\mathbb{N}}}

\def\Q{{\mathbb{Q}}}

\def\RP{{\mathbb{R}P^{n}}}

\def\aa{{\alpha}}
\def\bb{{\beta}}
\def\ga{{\gamma}}

\def\th{{\theta}}

\def\Lm{{\Lambda}}

\def\vf{{\varphi}}
\def\vs{{\varsigma}}

\def\ker{{\rm ker}}

\def\sgn{{\rm sgn}}

\def\hv{{\rm hv}}

\def\rank{{\rm rank}}

\def\mod{{\rm mod}}
\def\CG{{\rm CG}}

\def\dm{{\diamond}}
\def\ol#1{\overline{#1}}  
\def\td#1{\tilde{#1}}

\def\mapright#1{\smash{\mathop{\longrightarrow}\limits^{#1}}}

\def\mapdown#1{\Big\downarrow\rlap{$\vcenter{\hbox{$\scriptstyle#1$}}$}}

\def\mapright#1{\smash{\mathop{\longrightarrow}\limits^{#1}}}

\def\mapdown#1{\Big\downarrow\rlap{$\vcenter{\hbox{$\scriptstyle#1$}}$}}

\title{The existence of two non-contractible closed geodesics on every bumpy Finsler compact space form}
\author{Hui Liu$^{1}$,\thanks{Partially supported by NSFC (No. 11401555), Anhui Provincial Natural Science Foundation (No. 1608085QA01).
E-mail: huiliu00031514@whu.edu.cn. } \qquad  Yiming
Long$^{2, 3}$,\thanks{Partially supported by NSFC (Nos. 11131004 and 11671215), MCME and LPMC of
MOE of China, Nankai University and BAICIT of Capital Normal University. E-mail: longym@nankai.edu.cn.} \qquad
Yuming Xiao$^{4}$\thanks{Supported by the Scientific Research Funds for Young Teachers of Sichuan University, Grant 2012SCU11083.
e-mail: yumingxiao@scu.edu.cn.} \\ \\
$^{1}$ School of Mathematics and Statistics, Wuhan University,
\\Wuhan 430072, Hubei, China\\
$^{2}$ Chern Institute of Mathematics and LPMC, Nankai University, \\Tianjin 300071, China\\
$^{3}$ Beijing Advanced Innovation Center for Imaging Technology,\\Capital Normal University, Beijing 100048, China\\
$^{4}$ School of Mathematics, Sichuan University, Chengdu 610064, China \\}
\date{}
\begin{document}
\newtheorem{definition}{Definition}[section]
\newtheorem{theorem}{Theorem}[section]
\newtheorem{lemma}{Lemma}[section]
\newtheorem{corollary}{Corollary}[section]
\newtheorem{example}{Example}[section]
\newtheorem{property}{Property}[section]
\newtheorem{proposition}{Proposition}[section]
\newtheorem{remark}{Remark}[section]

\newcommand{\qed}{\nolinebreak\hfill\rule{2mm}{2mm}
\par\medbreak}
\newcommand{\Proof}{\par\medbreak\it Proof: \rm}
\newcommand{\rem}{\par\medbreak\it Remark: \rm}
\newcommand{\defi}{\par\medbreak\it Definition : \rm}
\renewcommand{\thefootnote}{\arabic{footnote}}

\maketitle

\begin{abstract}
{\it Let $M=S^n/ \Gamma$ and $h$ be a nontrivial element of finite order $p$ in $\pi_1(M)$, where the integer $n\geq2$,
$\Gamma$ is a finite group which acts freely and isometrically on the $n$-sphere and therefore $M$ is diffeomorphic to a compact space form. In this paper,
we establish first the resonance identity for non-contractible homologically visible minimal closed geodesics
of the class $[h]$ on every Finsler compact space form $(M, F)$ when there exist only finitely many distinct non-contractible closed
geodesics of the class $[h]$ on $(M, F)$. Then as an application of this resonance identity, we prove the existence of at least two distinct
non-contractible closed geodesics of the class $[h]$ on $(M, F)$ with a bumpy Finsler metric,
which improves a result of Taimanov in \cite{Tai2016} by removing some additional conditions. Also our results extend the resonance
identity and multiplicity results on $\RP$ in \cite{LX} to general compact space forms. }
\end{abstract}



\renewcommand{\theequation}{\thesection.\arabic{equation}}
\makeatletter\@addtoreset{equation}{section}\makeatother

\setcounter{equation}{0}
\section{Introduction}
Let $M=S^n/ \Gamma$ and $h$ be a nontrivial element of finite order $p$ in $\pi_1(M)$,
where the integer $n\geq2$, $\Gamma$ is a finite group which acts freely and
isometrically on the $n$-sphere and therefore $M$ is diffeomorphic to a compact space form
which is typically a non-simply connected manifold. In particular, if $\Gamma=\Z_2$,  then $S^n/ \Gamma$
is the $n$-dimensional real projective space $\mathbb{R}P^n$. Motivated by  the works \cite{XL2015},
\cite{DLX2015} and \cite{LX} about closed geodesics on Finsler $\mathbb{R}P^n$,
and based on Taimanov's work \cite{Tai2016} on rational equivariant cohomology of  non-contractible loops on $S^n/ \Gamma$,
this paper is concerned with the multiplicity of closed geodesics on Finsler $S^n/ \Gamma$.

Let $(M, F)$ be a Finsler manifold  and $\Lm M$ be the free loop space on $M$ defined by
\begin{equation*}
\label{LambdaM}
  \Lambda M=\left\{\gamma: S^{1}\to M\mid \gamma\ {\rm is\ absolutely\ continuous\ and}\
                        \int_{0}^{1}F(\gamma,\dot{\gamma})^{2}dt<+\infty\right\},
\end{equation*}
endowed with a natural structure of Riemannian Hilbert manifold on which the group $S^1=\R/\Z$ acts continuously by
isometries (cf. Shen \cite{Shen2001}). A closed geodesic $c:S^1=\R/\Z\to M$ is {\it prime}
if it is not a multiple covering (i.e., iteration) of any other closed geodesics. Here the $m$-th iteration $c^m$ of $c$ is defined
by $c^m(t)=c(mt)$.  The inverse curve $c^{-1}$ of $c$ is defined by
$c^{-1}(t)=c(1-t)$ for $t\in \R$.  Note that unlike Riemannian manifold,
the inverse curve $c^{-1}$ of a closed geodesic $c$
on a irreversible Finsler manifold need not be a geodesic.
We call two prime closed geodesics
$c$ and $d$ {\it distinct} if there is no $\th\in (0,1)$ such that
$c(t)=d(t+\th)$ for all $t\in\R$.
For a closed geodesic $c$ on $(M,\,F)$, denote by $P_c$
the linearized Poincar\'{e} map of $c$. Recall that a Finsler metric $F$ is {\it bumpy} if all the closed geodesics
on $(M, \,F)$ are non-degenerate, i.e., $1\notin \sigma(P_c)$ for any closed
geodesic $c$.

It is well known (cf. Chapter 1 of Klingenberg \cite{Kli1978}) that $c$ is a closed geodesic or a constant curve
on $(M,F)$ if and only if $c$ is a critical point of the energy functional
\begin{equation*}
\label{energy}
E(\gamma)=\frac{1}{2}\int_{0}^{1}F(\gamma,\dot{\gamma})^{2}dt.
\end{equation*}
Based on it, many important results on this subject have been obtained (cf. \cite{Ano},  \cite{Hin1984}-\cite{Hingston1993},
\cite{Rad1989}-\cite{Rad1992}). In particular, in 1969 Gromoll and Meyer \cite{GM1969JDG}
used Morse theory and Bott's index iteration formulae \cite{Bott1956}
to establish the existence of infinitely many distinct closed geodesics on $M$, when the Betti number sequence
$\{{\beta}_k(\Lm M;\mathbb{Q})\}_{k\in\Z}$ is unbounded. Then Vigu$\acute{e}$-Poirrier and Sullivan \cite{VS1976}
further proved in 1976 that for a compact simply connected manifold $M$, the Gromoll-Meyer condition holds if and
only if $H^{*}(M;\mathbb{Q})$ is generated by more than one element.

However, when $\{{\beta}_k(\Lm M;\mathbb{Q})\}_{k\in\Z}$ is bounded, the problem is quite complicated.
In 1973, Katok \cite{Katok1973} endowed some irreversible Finsler metrics to the compact rank one symmetric spaces
$$  S^{n},\ \mathbb{R}P^{n},\ \mathbb{C}P^{n},\ \mathbb{H}P^{n}\ \text{and}\ {\rm CaP}^{2},  \label{mflds}$$
each of which possesses only finitely many distinct prime closed geodesics (cf. also Ziller \cite{Ziller1977},\cite{Ziller1982}).
On the other hand, Franks \cite{Franks1992} and Bangert \cite{Bangert1993} together proved that
there are always infinitely many distinct closed geodesics on every Riemannian sphere $S^2$
(cf. also  Hingston \cite{Hingston1993}, Klingenberg \cite{Kli1995}).
These results imply that the metrics play an important role on the multiplicity of closed geodesics on those manifolds.

In 2004,  Bangert and Long \cite{BL2010} (published in 2010) proved the existence of at least two distinct
closed geodesics on every Finsler $S^2$. Subsequently, such a multiplicity result for
$S^{n}$ with a bumpy Finsler metric  was proved by Duan and Long
\cite{DL2007} and Rademacher \cite{Rad2010} independently. Furthermore  in a recent paper \cite{DLW1},
Duan, Long and Wang  proved the same conclusion for any compact simply-connected bumpy Finsler manifold.
We refer the readers to \cite{DL2010}-\cite{DLW2}, \cite{HiR}, \cite{LD2009}, \cite{Rad2007} \cite{Wang2008}-\cite{Wang2012}
and the references therein for more interesting results and the survey papers of Long \cite{lo2006},
Taimanov \cite{Tai2010}, Burns and Matveev \cite{BM2013} and Oancea \cite{Oancea2014} for more recent progresses on this subject.

Motivated by the studies on simply connected manifolds, in particular, the resonance identity proved by Rademacher \cite{Rad1989},
and based on Westerland's works \cite{West2005}, \cite{West2007} on loop homology of $\RP$, Xiao and Long \cite{XL2015}
in 2015 investigated the topological structure of the non-contractible
loop space and established the resonance identity for the non-contractible closed geodesics on $\R P^{2n+1}$ by use of
$\Z_2$ coefficient homology. As an application, Duan, Long and Xiao \cite{DLX2015} proved the existence of at least two
distinct non-contractible closed geodesics on $\R P^{3}$ endowed with a bumpy and irreversible Finsler metric. Subsequently in \cite{Tai2016}, Taimanov
used a quite different method from \cite{XL2015} to compute the rational equivariant cohomology of
 the non-contractible loop spaces in compact space forms $S^n/ \Gamma$ and
proved the existence of at least two distinct non-contractible closed geodesics on $\mathbb{R}P^2$ endowed with a bumpy and irreversible Finsler metric.
Then in \cite{Liu}, Liu combined Fadell-Rabinowitz index theory with Taimanov's topological results to get many multiplicity
results of non-contractible closed geodesics on positively curved Finsler $\RP$. Very recently,  Liu and Xiao \cite{LX} established
the resonance identity for the non-contractible closed geodesics on $\mathbb{R}P^n$, and
together with \cite{DLX2015} and \cite{Tai2016} proved the existence of at least two distinct
non-contractible closed geodesics on every bumpy $\mathbb{R}P^n$ with $n\geq2$.

Based on the works of \cite{DLW1} and \cite{LX}, it is natural to ask whether every bumpy Finsler compact space form
possesses two distinct closed geodesics on each of its nontrivial classes. This paper gives a positive answer to this question.
To this end, we first establish the following resonance identity in section 2.
Comparing with Theorem 1.1 of \cite{LX}, the difficulties mainly lie in that the parity of the order $p$ of the nontrivial element $h$ in $\pi_1(M)$
is unknown which yields that the computations of critical modules of non-contractible closed geodesics are
very complicated (cf. Lemma \ref{Rad1992} below) and the parity of $i(c^{p+1})-i(c)$ is unknown for any closed geodesic $c$,
and the proofs of the positivity of mean index of non-contractible homologically visible minimal closed
geodesics for compact space forms(cf. Lemma \ref{Lm2.2} below), and also a non-contractible minimal closed geodesic $c$ of the class $[h]$
may be some iteration of a closed geodesic $\gamma$ which is not in the class $[h]$. Recall that $\Gamma$ is a finite group which acts freely and
isometrically on the $n$-sphere.
\begin{theorem}\label{Thm1.1} Let $M=S^n/ \Gamma$ and $h$ be a nontrivial element of finite order $p$ in $\pi_1(M)$.
Suppose the Finsler manifold $(M,F)$ possesses only finitely
many distinct non-contractible minimal closed geodesics of the class $[h]$,
among which we denote the distinct non-contractible homologically visible minimal
closed geodesics by $c_1, \ldots, c_r$  for some integer $r>0$, where $n\geq2$ and
a closed geodesic $c$ of the class $[h]$ is called {\it minimal} if it is not an iteration of any other
closed geodesics in class $[h]$. Then we have
\bea  \sum_{j=1}^{r}\frac{\hat{\chi}(c_j)}{\hat{i}(c_j)} = \bar{B}(\Lm_hM ;\Q) =\left\{\begin{array}{ll}
    \frac{n+1}{2(n-1)},&\ if\ n\in 2\N-1,\\
    \frac{n}{2(n-1)},&\ if\ n\in 2\N.\\
    \end{array}\right.     \lb{reident1}\eea
where the mean Euler number $\hat{\chi}(c_j)$ of $c_j$ is defined by
$$  \hat{\chi}(c_j) = \frac{1}{n_j}\sum_{m=1}^{n_j/p}\sum_{l=0}^{2n-2}(-1)^{l+i(c_{j}^{p(m-1)+1})}k_{l}^{\epsilon(c_j^{p(m-1)+1})}(c_{j}^{p(m-1)+1})\in\Q, $$
and $n_j=n_{c_j}$ is the analytical period of $c_j$, $k_{l}^{\epsilon(c_j^{p(m-1)+1})}(c_{j}^{p(m-1)+1})$
is the local homological type number of $c_{j}^{p(m-1)+1}$,
$i(c_{j})$ and $\hat{i}(c_j)$ are the Morse index and mean index of $c_j$ respectively.

In particular, if the Finsler metric $F$ on $M=S^n/ \Gamma$ is bumpy, then (\ref{reident1}) has the following simple form
\bea  \sum_{j=1}^{r}\left((-1)^{i(c_{j})}k_0^{\epsilon(c_j)}(c_{j})+(-1)^{i(c_{j}^{p+1})}k_0^{\epsilon(c_j^{p+1})}(c_{j}^{p+1})\right)\frac{1}{\hat{i}(c_{j})}
=\left\{\begin{array}{ll}
    \frac{p(n+1)}{n-1},&\ if\ n\in 2\N-1,\\
    \frac{pn}{n-1},&\ if\ n\in 2\N.\\
    \end{array}\right.  \lb{breident1}\eea
\end{theorem}

Based on Theorem \ref{Thm1.1},  we use Morse theory and draw support from the well known Kronecker approximation theorem
to prove our main multiplicity result of non-contractible closed geodesics on $(S^{2n+1}/ \Gamma, F)$.

\begin{theorem}
\label{mainresult}
Let $M=S^{2n+1}/ \Gamma$ and $h$ be a nontrivial element of finite order $p$ in $\pi_1(M)$. Then every bumpy Finsler metric $F$ on $M$
has at least two distinct non-contractible closed geodesics of the class $[h]$.
\end{theorem}

Note that the only non-trivial group which acts freely on $S^{2n}$ is $\Z_2$ and $S^{2n}/ \Z_2=\mathbb{R}P^{2n}$(cf. P.5 of \cite{Tai2016}). Since we
have proved the same result as the above Theorem 1.2 for $\mathbb{R}P^{2n}$ in Theorem 1.2 and Corollary 1.1 of \cite{LX}, then we have

\begin{theorem}
\label{mainresult2}
Let $M=S^{n}/ \Gamma$ and $h$ be a nontrivial element of finite order $p$ in $\pi_1(M)$, where $n\geq2$. Then every bumpy Finsler metric $F$ on $M$
has at least two distinct non-contractible closed geodesics of the class $[h]$.
\end{theorem}

\begin{remark}

(i) In Theorem 5 of \cite{Tai2016}, Taimanov
proved the same result as Theorem 1.2 under the conditions that $\pi_1(\Lm_{h} (M))_{SO(2)}\neq 1$,
$h$ has an even order in $\pi_1(M)$ and the centralizer of $h$ are pairwise non-conjugate,
our Theorem 1.2 improves Taimanov's result by removing these additional conditions.

(ii) When $\Gamma=\Z_2$, then $S^n/ \Gamma$ is the $n$-dimensional real
projective space $\mathbb{R}P^n$ and $p=2$, one can easily check that for $\mathbb{R}P^n$,
the results of the above Theorems 1.1-1.3 are just the results of Theorems 1.1-1.2
and Corollary 1.1 of \cite{LX}. So the main results of this paper are generalizations of those of \cite{LX}.
Note that the only non-trivial group which acts freely on $S^{2n}$ is $\Z_2$ and
$S^{2n}/ \Z_2=\mathbb{R}P^{2n}$(cf. P.5 of \cite{Tai2016}), then we only need to prove
Theorem 1.1 for the case when $n$ is odd.

 \end{remark}

This paper is organized as follows. In section 2, we apply Morse theory to the non-contractible loops of the class $[h]$
and  establish the resonance identity of Theorem 1.1.
Then in section 3, we firstly recall the precise iteration formulae of Morse indices for orientable closed geodesics,
and combine it with Theorem 1.1 to investigate the Morse indices for closed geodesics on $S^n/ \Gamma$
and build a bridge between the Morse indices and a division of an interval, then our problem are reduced to
a problem in Number Theory and we review some theories about a special system of irrational numbers associated to our problem
developed in \cite{LX}.
In section 4, we draw support from the well known Kronecker's
approximation theorem and other techniques in Number theory to give the proof of Theorem \ref{mainresult}.
Finally in section 5, for the reader's convenience, we give the proof of Theorem \ref{crucialtheorem}
about a special system of irrational numbers as an appendix.

In this paper, let $\N$, $\N_0$,  $\Z$, $\Q$ and $\Q^{c}$ denote the sets of natural integers, non-negative
integers, integers, rational numbers and  irrational numbers respectively. We also use  notations $E(a)=\min\{k\in\Z\,|\,k\ge a\}$,
$[a]=\max\{k\in\Z\,|\,k\le a\}$, $\varphi(a)=E(a)-[a]$ and $\{a\}=a-[a]$ for any $a\in\R$.
Throughout this paper, we use $\Q$ coefficients for all homological and cohomological modules.

\setcounter{equation}{0}
\section{Resonance identity of non-contractible closed geodesics on $(S^n/ \Gamma,F)$}
Let $M=(M,F)$ be a compact Finsler manifold, the space
$\Lambda=\Lambda M$ of $H^1$-maps $\gamma:S^1\rightarrow M$ has a
natural structure of Riemannian Hilbert manifolds on which the
group $S^1=\R/\Z$ acts continuously by isometries. This action is defined by
$(s\cdot\gamma)(t)=\gamma(t+s)$ for all $\gamma\in\Lm$ and $s,
t\in S^1$. For any $\gamma\in\Lambda$, the energy functional is
defined by
\be E(\gamma)=\frac{1}{2}\int_{S^1}F(\gamma(t),\dot{\gamma}(t))^2dt.
\lb{2.1}\ee
It is $C^{1,1}$ and invariant under the $S^1$-action. The
critical points of $E$ of positive energies are precisely the closed geodesics
$\gamma:S^1\to M$. The index form of the functional $E$ is well
defined along any closed geodesic $c$ on $M$, which we denote by
$E''(c)$. As usual, we denote by $i(c)$ and
$\nu(c)$ the Morse index and nullity of $E$ at $c$. In the
following, we denote by
\be \Lm^\kappa=\{d\in \Lm\;|\;E(d)\le\kappa\},\quad \Lm^{\kappa-}=\{d\in \Lm\;|\; E(d)<\kappa\},
  \quad \forall \kappa\ge 0. \lb{2.2}\ee
For a closed geodesic $c$ we set $ \Lm(c)=\{\ga\in\Lm\mid E(\ga)<E(c)\}$.

For $m\in\N$ we denote the $m$-fold iteration map
$\phi_m:\Lambda\rightarrow\Lambda$ by $\phi_m(\ga)(t)=\ga(mt)$, for all
$\,\ga\in\Lm, t\in S^1$, as well as $\ga^m=\phi_m(\gamma)$. If $\gamma\in\Lambda$
is not constant then the multiplicity $m(\gamma)$ of $\gamma$ is the order of the
isotropy group $\{s\in S^1\mid s\cdot\gamma=\gamma\}$. For a closed geodesic $c$,
the mean index $\hat{i}(c)$ is defined as usual by
$\hat{i}(c)=\lim_{m\to\infty}i(c^m)/m$. Using singular homology with rational
coefficients we consider the following critical $\Q$-module of a closed geodesic
$c\in\Lambda$:
\be \overline{C}_*(E,c)
   = H_*\left((\Lm(c)\cup S^1\cdot c)/S^1,\Lm(c)/S^1; \Q\right). \lb{2.3}\ee
In the following we let $M=S^n/ \Gamma$ and $h$ be a nontrivial element of finite order $p$ in $\pi_1(M)$, where the integer $n\geq2$,
$\Gamma$ acts freely and isometrically on the $n$-sphere and therefore
$M$ is diffeomorphic to a compact space form. Then the free loop space $\Lambda M$ possesses
a natural decomposition\bea \Lambda M=\bigsqcup_{g\in \pi_1(M)}\Lambda_g M,\nn\eea
where $\Lambda_g M$ is the connected components of $\Lambda M$ whose elements are homotopic
to $g$. We set $\Lambda_{h}(c) = \{\gamma\in \Lambda_{h}M\mid E(\gamma)<E(c)\}$.
Note that for a non-contractible minimal closed geodesic $c$ of class $[h]$, $c^m\in\Lambda_h M$ if and only if $m\equiv 1(\mod~ p)$.

We call a non-contractible minimal closed geodesic $c$ of class $[h]$ satisfying the isolation condition, if
the following holds:

{\bf (Iso)  For all $m\in\N$ the orbit $S^1\cdot c^{p(m-1)+1}$ is an
isolated critical orbit of $E$. }

Note that if the number of non-contractible minimal closed geodesics of class $[h]$ on $M$
is finite, then all the non-contractible minimal closed geodesics of class $[h]$ satisfy (Iso).

For a non-contractible closed geodesic $d$ of class $[h]$, we
have $d=c^{p(m-1)+1}$ for some $m\in\N$, where $c$ is a minimal closed geodesic of class $[h]$
and $c=\gamma^t$ for a prime closed geodesic $\gamma$ with $t\in\N$ .
Then $d$ has multiplicity $tp(m-1)+t$, the subgroup $\Z_{tp(m-1)+t}=\{\frac{l}{tp(m-1)+t}\mid 0\leq l<tp(m-1)+t\}$
of $S^1$ acts on $\overline{C}_*(E, d)$. As studied in p.59 of \cite{Rad1992},
for all $m\in\N$, let
$H_{\ast}(X,A)^{\pm\Z_{tp(m-1)+t}}
   = \{[\xi]\in H_{\ast}(X,A)\,|\,T_{\ast}[\xi]=\pm [\xi]\}$,
where $T$ is a generator of the $\Z_{tp(m-1)+t}$-action.
On $S^1$-critical modules of $c^{p(m-1)+1}$, the following lemma holds:
\begin{lemma}
\label{Rad1992} {\rm (cf. Satz 6.11 of \cite{Rad1992} and \cite{BL2010})}  Suppose $c$ is
a non-contractible minimal closed geodesic of class $[h]$ on a Finsler manifold $M=S^n/ \Gamma$ satisfying (Iso). Then
there exist $U_{\gamma^{tp(m-1)+t}}$ and $N_{\gamma^{tp(m-1)+t}}$, the so-called local negative
disk and the local characteristic manifold at $c^{p(m-1)+1}$ respectively,
such that $\nu(c^{p(m-1)+1})=\dim N_{\gamma^{tp(m-1)+t}}$ and
\bea &&\overline{C}_q( E,c^{p(m-1)+1})
\equiv H_q\left((\Lm_h(c^{p(m-1)+1})\cup S^1\cdot c^{p(m-1)+1})/S^1, \Lm_h(c^{p(m-1)+1})/S^1\right)\nn\\
&=& (H_{i(c^{p(m-1)+1})}(U_{\gamma^{tp(m-1)+t}}^-\cup\{\gamma^{tp(m-1)+t}\},U_{\gamma^{tp(m-1)+t}}^-)\nn\\
    &&\otimes H_{q-i(c^{p(m-1)+1})}(N_{\gamma^{tp(m-1)+t}}^-\cup\{\gamma^{tp(m-1)+t}\},N_{\gamma^{tp(m-1)+t}}^-))^{+\Z_{tp(m-1)+t}}, \nn
\eea

(i) When $\nu(c^{p(m-1)+1})=0$, there holds
\bea \overline{C}_q( E,c^{p(m-1)+1}) = \left\{\begin{array}{ll}

     \Q, &\quad {\it if}\;i(c^{p(m-1)+1})-i(\gamma)\in 2\Z\;and\;
                   q=i(c^{p(m-1)+1}),\;  \cr
     0, &\quad {\it otherwise},\\ \end{array}\right. \nn \eea

(ii) When $\nu(c^{p(m-1)+1})>0$, there holds
$$ \overline{C}_q( E,c^{p(m-1)+1})=H_{q-i(c^{p(m-1)+1})}(N_{\gamma^{tp(m-1)+t}}^-\cup\{\gamma^{tp(m-1)+t}\},
N_{\gamma^{tp(m-1)+t}}^-)^{\epsilon(c^{p(m-1)+1})\Z_{tp(m-1)+t}}, $$
where $\epsilon(c^{p(m-1)+1})=(-1)^{i(c^{p(m-1)+1})-i(\gamma)}$.
\end{lemma}

As usual, for $m\in\N$ and $l\in\Z$ we define the local homological type numbers of $c^{p(m-1)+1}$ by
\be k_{l}^{\epsilon(c^{p(m-1)+1})}(c^{p(m-1)+1})
= \dim H_{l}(N_{\gamma^{tp(m-1)+t}}^-\cup\{\gamma^{tp(m-1)+t}\},N_{\gamma^{tp(m-1)+t}}^-)^{\epsilon(c^{p(m-1)+1})\Z_{tp(m-1)+t}}.  \lb{CGht1}\ee

Based on works of Rademacher in \cite{Rad1989}, Long and Duan in \cite{LD2009} and \cite{DL2010},
we define the {\it analytical period} $n_c$ of the closed geodesic $c$ by
\be n_c = \min\{j\in 2p\N\,|\,\nu(c^j)=\max_{m\ge 1}\nu(c^m)
                  \}. \lb{CGap1}\ee
Note that here in order to simplify the study for non-contractible closed geodesics of class $[h]$ on $M=S^n/ \Gamma$,
we have slightly modified the definition in \cite{LD2009} and \cite{DL2010} by requiring the analytical
period to be integral multiple of $2p$. Then by the same proofs in \cite{LD2009} and \cite{DL2010}, we have
\be  k_{l}^{\epsilon(c^{p(m-1)+1+kn_c})}(c^{p(m-1)+1+kn_c})
= k_{l}^{\epsilon(c^{p(m-1)+1})}(c^{p(m-1)+1}), \qquad \forall\;m,\;k\in \N,\;l\in\Z.  \lb{CGap2}\ee
For more detailed properties of the analytical period $n_c$ of a closed geodesic $c$, we refer readers to
the two Section 3s in \cite{LD2009} and \cite{DL2010}.

As in \cite{BK1983}, we have
\begin{definition}\label{def-hv} Let $(M,F)$ be a compact Finsler manifold. A closed geodesic $c$ on $M$
is homologically visible, if there exists an integer $k\in\Z$ such that $\bar{C}_k(E,c) \not= 0$.
\end{definition}
\begin{lemma}\label{Lm2.2} Suppose the Finsler manifold $M=S^n/ \Gamma$ possesses only finitely
many distinct non-contractible minimal closed geodesics of the class $[h]$, among
which we denote the distinct non-contractible homologically visible minimal closed
geodesics by $c_1, \ldots, c_r$  for some integer $r>0$. Then we have
\be  \hat{i}(c_i) > 0, \qquad \forall\;1\le i\le r.  \lb{fcg.1}\ee
\end{lemma}

\Proof First, we claim that Theorem 3 in \cite{BK1983} for $M=S^n/ \Gamma$ can be stated as:

`` Let $c$ be a closed geodesic in $\Lambda_{h}M$ such that $i(c^{m})=0$ for all $m\in\N$. Suppose $c$ is neither
homologically invisible nor an absolute minimum of $E$ in $\Lambda_{h}M$. Then there exist infinitely many closed geodesics in $\Lambda_{h}M$."

Indeed, one can focus the proofs of Theorem 3 in \cite{BK1983} on $\Lambda_{h}M$ with some obvious modifications.
Assume by contradiction. Similarly as in \cite{BK1983}, we can choose a different $c\in\Lm_{h}M$, if necessary, and find $p\in\N$ such that
$H_{p}(\Lambda_{h}(c)\cup S\cdot c,\Lambda_{h}(c))\neq0$ and $H_{q}(\Lambda_{h}(c)\cup S\cdot c,\Lambda_{h}(c))=0$ for every $q>p$ and
every closed geodesic $d\in\Lm_{h}M$ with $i(d^{m})\equiv0.$

Consider the following commutative diagram \be
\begin{tabular}{ccc}
$H_{p}(\Lambda_{h}(c)\cup S\cdot c,\Lambda_{h}(c))$   & $\mapright{{\psi}^{m}_{*}}$ & $H_{p}(\Lambda_{h}(c^{m})\cup S\cdot c^{m},\Lambda_{h}(c^{m}))$ \\
$\mapdown{i_*}$&                       & $\mapdown{i_*}$ \\
$H_{p}(\Lambda_{h}M,\Lambda_{h}(c))$                     & $\mapright{{\psi}^{m}_{*}}$ & $H_{p}(\Lambda_{h}M,\Lambda_{h}(c^{m}))$,\\
\end{tabular}  \lb{926diagram}\ee
where $m\equiv 1(\mod~ p)$ and $\psi^{m}:\Lm_{h} M\to\Lm_{h} M$ is the $m$-fold iteration map. By similar arguments as those in \cite{BK1983},
there is $A>0$ such that the map $i_{*}\circ\psi^{m}_{*}$ is one-to-one, if $E(c^{m})>A$ and none of the $k_{i}\in K_0$ divides $m$
where  $$K_{0}=\{k_{0},k_{1},k_{2},\dots,k_{s}\},$$ with $k_{0}=p$ and $k_{1},k_{2},\dots,k_{s}$ therein. Here note that the
required $m\equiv 1(\mod~ p)$ and so $c^{m}\in\Lambda_{h}(M)$ for $c\in\Lambda_{h}M$.

On the other hand, we define $$K=\{m\geq2\mid E(c^{m})\leq A\}\cup K_{0}.$$
Then by Corollary 1 of \cite{BK1983}, there exists $\bar{m}\in\N\backslash\{1\}$ such that no $k\in K$ divides $\bar{m}$ and
$\psi^{\bar{m}}_{*}\circ i_{*}$ vanishes.  In particular, $E(c^{\bar{m}})>A$ and none of the $k_{i}\in K_0$ divides $\bar{m}$.
Due to $\psi^{\bar{m}}_{*}\circ i_{*}=i_{*}\circ \psi^{\bar{m}}_{*}$ in (\ref{926diagram}),  this yields a contradiction.
Hence there exist infinitely many closed geodesics in $\Lambda_{h}M$.

Accordingly, Corollary 2 in \cite{BK1983}  for $M=S^n/ \Gamma$ can be stated as:

``  Suppose there exists a closed geodesic $c\in\Lambda_{h}M$ such that $c^{m}$ is a local minimum of $E$ in $\Lambda_{h}M$ for
infinitely many $m\equiv 1(\mod~ p)$. Then there exist infinitely many closed geodesics in $\Lambda_{h}M$."

Based on the above two variants of Theorem 3 and Corollary 2 in \cite{BK1983}, we can prove our Lemma \ref{Lm2.2} as follows.

It is well known that every closed geodesic $c$ on $M$ must have mean index $\hat{i}(c)\ge 0$.

Assume by contradiction that there is a non-contractible homologically visible minimal closed geodesic $c$ of the class $[h]$ on $M$ satisfying
$\hat{i}(c)=0$. Then $i(c^m)=0$ for all $m\in\N$ by Bott iteration formula and $c$ must be an absolute
minimum of $E$ in $\Lambda_{h}M$, since otherwise there would exist infinitely many distinct non-contractible
closed geodesics of the class $[h]$ on $M$ by the above variant of Theorem 3 on p.385 of \cite{BK1983}.

On the other hand, by Lemma 7.1 of \cite{Rad1992}, there exists a $k(c)\in p\N$ such that
$\nu(c^{m+k(c)})=\nu(c^m)$ for all $m\in\N$. Specially we obtain $\nu(c^{mk(c)+1}) = \nu(c)$ for all
$m\in\N$ and then elements of $\ker(E''(c^{mk(c)+1}))$ are precisely $mk(c)+1$st iterates of
elements of $\ker(E''(c))$. Thus by the Gromoll-Meyer theorem in \cite{GM1969Top}, the behavior of the
restriction of $E$ to $\ker(E''(c^{mk(c)+1}))$ is the same as that of the restriction of $E$ to
$\ker(E''(c))$. Then together with the fact $i(c^m)=0$ for all $m\in\N$, we obtain that $c^{mk(c)+1}$
is a local minimum of $E$ in $\Lambda_{h}M$ for every $m\in\N$. Because $M$ is compact and
possessing finite fundamental group ($\pi_1(M)$ is finite for the spherical space forms!), there must exist infinitely
many distinct non-contractible closed geodesics of the class $[h]$ on $M$
by the above variant of Corollary 2 on p.386 of \cite{BK1983}. Then it yields a contradiction and proves (\ref{fcg.1}).
$\hfill\Box$

In \cite{Tai2016}, Taimanov calculated the rational equivariant cohomology of the spaces of non-contractible loops in compact space forms
which is crucial for us to prove Theorem \ref{Thm1.1} and  can be stated as follows.
\begin{lemma}\label{Lm2.3} {\rm (cf. Theorem 3 of \cite{Tai2016})} For $M=S^n/ \Gamma$, we have

(i) When $n=2k+1$ is odd, the $S^1$-cohomology ring of $\Lambda_h M$ has the form
$$H^{S^1, *}(\Lambda_h M; \Q)=\Q[w, z]/ \{w^{k+1} = 0\}, \quad deg(w)=2, \quad deg(z)=2k$$
Then the $S^1$-equivariant Poincar$\acute{e}$ series
of $\Lambda_h M$ is given by
\bea P^{S^1}(\Lambda_h M; \Q)(t)&=&\frac{1-t^{2k+2}}{(1-t^2)(1-t^{2k})}\nn\\
&=&\frac{1}{1-t^2}+\frac{t^{2k}}{1-t^{2k}}\nn\\&=&(1+t^2+t^4+\cdots+t^{2k}+\cdots)+(t^{2k}+t^{4k}+t^{6k}+\cdots),\nn\eea
which yields Betti numbers
\bea \bar{\beta}_q=\rank H_q^{S^1}(\Lambda_h M;\Q)=\left\{\begin{array}{ll}
    2,&\ if\ q\in \{j(n-1)\mid j\in\mathbb{N}\},\\
    1,&\ if\ q\in (2\mathbb{N}_0)\backslash\{j(n-1)\mid j\in\mathbb{N}\},\\
    0,&\ otherwise.\\
    \end{array}\right. \lb{b.1}\eea
and the average $S^1$-equivariant Betti number of $\Lambda_h M$ satisfies
\be  \bar{B}(\Lambda_h M;\Q)\equiv\lim_{q\to+\infty}\frac{1}{q}\sum_{k=0}^{q}(-1)^k\bar{\beta}_k= \frac{n+1}{2(n-1)}.   \lb{aB.1}\ee

(ii) When $n=2k$ is even, the $S^1$-cohomology ring of $\Lambda_h M$ has the form
$$H^{S^1, *}(\Lambda_h M; \Q)=\Q[w, z]/ \{w^{2k} = 0\}, \quad deg(w)=2, \quad deg(z)=4k-2$$
Then the $S^1$-equivariant Poincar$\acute{e}$ series
of $\Lambda_h M$ is given by \bea P^{S^1}(\Lambda_h M; \Q)(t)&=&\frac{1-t^{4k}}{(1-t^2)(1-t^{4k-2})}\nn\\
&=&\frac{1}{1-t^2}+\frac{t^{4k-2}}{1-t^{4k-2}}\nn\\&=&(1+t^2+t^4+\cdots+t^{2k}+\cdots)+(t^{4k-2}+t^{2(4k-2)}+t^{3(4k-2)}+\cdots),\nn\eea
which yields Betti numbers \bea \bar{\beta}_q= \rank H_q^{S^1}(\Lambda_h M;\Q)=\left\{\begin{array}{ll}
2,&\ if\ q\in \{2j(n-1)\mid j\in\mathbb{N}\},\\
1,&\ if\ q\in (2\mathbb{N}_0)\backslash\{2j(n-1)\mid j\in\mathbb{N}\},\\
0,&\ otherwise.\\
\end{array}
\right.\lb{b.2}\eea
and the average $S^1$-equivariant Betti number of $\Lambda_h M$ satisfies
\be  \bar{B}(\Lambda_h M;\Q)\equiv\lim_{q\to+\infty}\frac{1}{q}\sum_{k=0}^{q}(-1)^k\bar{\beta}_k= \frac{n}{2(n-1)}.   \lb{aB.2}\ee
\end{lemma}

Now we give the proof of the resonance identity in Theorem \ref{Thm1.1}.

{\bf Proof of Theorem \ref{Thm1.1}.} Recall that we denote the non-contractible homologically visible minimal closed
geodesics of the class $[h]$ by $\CG_{\hv}^{[h]}(M)=\{c_1, \ldots, c_r\}$
for some integer $r>0$ when the number of distinct non-contractible minimal
closed geodesics of the class $[h]$ on $M=S^n/ \Gamma$ is finite. Note also that by Lemma \ref{Lm2.2} we have $\hat{i}(c_j)>0$ for all
$1\le j\le r$. In the following proof of Theorem \ref{Thm1.1}, we assume $n=2k+1$ for $k\in\N$ by Remark 1.1 (iii), then $M$ is orientable.

Let \bea m_q\equiv M_q(\Lm_hM) =\sum_{1\le j\le r,\; m\ge 1}\dim{\ol{C}}_q(E, c^{p(m-1)+1}_j),\quad q\in\Z.\nn\eea
The Morse series of $\Lm_hM$ is defined by
\be  M(t) = \sum_{q=0}^{+\infty}m_qt^q.  \label{wh}\ee

{\bf Claim 1.} {\it $\{m_q\}$ is a bounded sequence.}

In fact, by (\ref{CGap2}), we have \bea m_q=\sum_{j=1}^{r}\sum_{m=1}^{n_j/p}\sum_{l=0}^{2n-2}k_{l}^{\epsilon(c_j^{p(m-1)+1})}(c_{j}^{p(m-1)+1})
            \;{}^{\#}\left\{s\in\mathbb{N}_0\mid q-i(c_{j}^{p(m-1)+1+sn_j})=l\right\}, \label{wh2}\eea
by Theorem 10.1.2 of \cite{lo2000} and Lemma \ref{orient} below,
we have $|i(c_{j}^{p(m-1)+1+sn_j})-(p(m-1)+1+sn_j)\hat{i}(c_{j})|\le n-1$, then
\bea
&&{}^{\#}\left\{s\in\mathbb{N}_0\mid q-i(c_{j}^{p(m-1)+1+sn_j})=l\right\} \nn\\
&&\qquad = \;{}^{\#}\left\{s\in\mathbb{N}_0\mid l+i(c_{j}^{p(m-1)+1+sn_j})=q,\;
           |i(c_{j}^{p(m-1)+1+sn_j})-(p(m-1)+1+sn_j)\hat{i}(c_{j})|\le n-1\right\}  \nn\\
&&\qquad \le \;{}^{\#}\left\{s\in\mathbb{N}_0\mid n-1\geq |q-l-(p(m-1)+1+sn_{j})\hat{i}(c_{j})|  \right\}  \nn\\
&&\qquad = \;{}^{\#}\left\{s\in\mathbb{N}_0\mid \frac{q-l-n+1-(p(m-1)+1)\hat{i}(c_j)}{n_j\hat{i}(c_j)} \leq s
       \le \frac{q-l+n-1-(p(m-1)+1)\hat{i}(c_j)}{n_j\hat{i}(c_j)}\right\} \nn\\
&&\qquad \le \frac{2n-2}{n_j\hat{i}(c_j)} + 1.  \label{wh3}\eea
Hence Claim 1 follows by (\ref{wh2}) and (\ref{wh3}).

We now use the method
in the proof of Theorem 5.4 of \cite{LW2007} to estimate
$$   M^{q}(-1) = \sum_{k=0}^{q}m_k(-1)^k.   $$

By (\ref{wh2}), Lemma \ref{orient} below and the fact that $n_j\in2\N$, we obtain
\bea M^{q}(-1)
&=& \sum_{k=0}^{q}m_{k}(-1)^{k}  \nn\\
&=& \sum_{j=1}^{r}\sum_{m=1}^{n_j/p}\sum_{l=0}^{2n-2}\sum_{k=0}^{q}(-1)^{k}k_{l}^{\epsilon(c_j^{p(m-1)+1})}(c_{j}^{p(m-1)+1})
            \;{}^{\#}\left\{s\in\mathbb{N}_0\mid k-i(c_{j}^{p(m-1)+1+sn_j})=l\right\}  \nn\\
&=& \sum_{j=1}^{r}\sum_{m=1}^{n_j/p}\sum_{l=0}^{2n-2}(-1)^{l+i(c_{j}^{p(m-1)+1})}k_{l}^{\epsilon(c_j^{p(m-1)+1})}(c_{j}^{p(m-1)+1})
            \;{}^{\#}\left\{s\in\mathbb{N}_0\mid l+i(c_{j}^{p(m-1)+1+sn_j})\le q\right\}.  \nn\eea
On the one hand, we have
\bea
&&{}^{\#}\left\{s\in\mathbb{N}_0\mid l+i(c_{j}^{p(m-1)+1+sn_j})\le q\right\} \nn\\
&&\qquad = \;{}^{\#}\left\{s\in\mathbb{N}_0\mid l+i(c_{j}^{p(m-1)+1+sn_j})\le q,\;
           |i(c_{j}^{p(m-1)+1+sn_j})-(p(m-1)+1+sn_j)\hat{i}(c_{j})|\le n-1\right\}  \nn\\
&&\qquad \le \;{}^{\#}\left\{s\in\mathbb{N}_0\mid 0\le (p(m-1)+1+sn_{j})\hat{i}(c_{j})\le q-l+n-1  \right\}  \nn\\
&&\qquad = \;{}^{\#}\left\{s\in\mathbb{N}_0\mid 0 \leq s
       \le \frac{q-l+n-1-(p(m-1)+1)\hat{i}(c_j)}{n_j\hat{i}(c_j)}\right\} \nn\\
&&\qquad \le \frac{q-l+n-1}{n_j\hat{i}(c_j)} + 1. \nn\eea
On the other hand, we have
\bea
&&{}^{\#}\left\{s\in\mathbb{N}_0\mid l+i(c_{j}^{p(m-1)+1+sn_j})\le q \right\} \nn\\
&&\qquad = \;{}^{\#}\left\{s\in\mathbb{N}_0\mid l+i(c_{j}^{p(m-1)+1+sn_j})\le q,\;
           |i(c_{j}^{p(m-1)+1+sn_j})-(p(m-1)+1+sn_{j})\hat{i}(c_{j})|\le n-1\right\}  \nn\\
&&\qquad \ge \;{}^{\#}\left\{s\in\mathbb{N}_0\mid i(c_{j}^{p(m-1)+1+sn_j})
       \le (p(m-1)+1+sn_{j})\hat{i}(c_{j})+n-1\le q-l \right\} \nn\\
&&\qquad \ge \;{}^{\#}\left\{s\in\mathbb{N}_0\mid 0 \le s
       \le \frac{q-l-n+1-(p(m-1)+1)\hat{i}(c_{j})}{n_j\hat{i}(c_{j})} \right\}  \nn\\
&&\qquad \ge \frac{q-l-n+1}{n_j\hat{i}(c_{j})} - 1.  \nn\eea
Thus we obtain
$$  \lim_{q\to+\infty}\frac{1}{q}M^{q}(-1)
  = \sum_{j=1}^{r}\sum_{m=1}^{n_j/p}\sum_{l=0}^{2n-2}(-1)^{l+i(c_{j}^{p(m-1)+1})}k_{l}^{\epsilon(c_j^{p(m-1)+1})}(c_{j}^{p(m-1)+1})\frac{1}{n_{j}\hat{i}(c_{j})}
              = \sum_{j=1}^{r}\frac{\hat{\chi}(c_j)}{\hat{i}(c_j)}.  $$
Since $m_{q}$ is bounded, we then obtain
$$  \lim_{q\to+\infty}\frac{1}{q}M^{q}(-1) = \lim_{q\to+\infty}\frac{1}{q}P^{S^{1},q}(\Lambda_{h}M;\Q)(-1)
          = \lim_{q\to+\infty}\frac{1}{q}\sum_{k=0}^{q}(-1)^k\bar{\beta}_k = \bar{B}(\Lm_hM;\Q),  $$
where $P^{S^{1},q}(\Lambda_{h}M;\Q)(t)$ is the truncated polynomial of
$P^{S^{1}}(\Lambda_{h}M;\Q)(t)$ with terms of degree less than or equal to $q$.
Thus by (\ref{aB.1}) we get
$$  \sum_{j=1}^{r}\frac{\hat{\chi}(c_j)}{\hat{i}(c_j)} =
    \frac{n+1}{2(n-1)},\forall\ n\in 2\N+1. $$
which proves (\ref{reident1}) of Theorem \ref{Thm1.1}. For the special case when each $c_{j}^{p(m-1)+1}$
is non-degenerate with $1\leq j\leq r$ and $m\in\mathbb{N}$,
we have $n_{j}=2p$ and $k_{l}^{\epsilon(c_j^{p(m-1)+1})}(c_{j}^{p(m-1)+1})=1$ when $l=0$ and $i(c_j^{p(m-1)+1})-i(\gamma_j)\in 2\Z$,
where $c_j$ is some iteration of a prime closed geodesic $\gamma_j$,
and $k_{l}^{\epsilon(c_j^{p(m-1)+1})}(c_{j}^{p(m-1)+1})=0$ for all other
$l\in\Z$. Then (\ref{reident1}) has the following simple form
\be  \sum_{j=1}^{r}\left((-1)^{i(c_{j})}k_0^{\epsilon(c_j)}(c_{j})+(-1)^{i(c_{j}^{p+1})}k_0^{\epsilon(c_j^{p+1})}(c_{j}^{p+1})\right)\frac{1}{\hat{i}(c_{j})}=
    \frac{p(n+1)}{n-1},\forall\ n\in 2\N+1.  \nn\ee
which proves (\ref{breident1}) of Theorem \ref{Thm1.1}. $\hfill\Box$

\setcounter{equation}{0}
\section{Preliminary for the proof of Theorem 1.2}
\subsection{Index iteration formulae for closed geodesics}
In \cite{lo1999} of 1999, Y. Long established the basic normal form
decomposition of symplectic matrices. Based on it, he
further established the precise iteration formulae of Maslov $\omega$-indices for
symplectic paths in \cite{lo2000}, which can be related to  Morse indices
of either orientable or non-orientable closed geodesics in a slightly different way (cf. \cite{Liu2005}, \cite{LL2002}
and Chap. 12 of \cite{lo2002}). Roughly speaking, the orientable (resp. non-orientable) case  corresponds to $i_{1}$ (resp. $i_{-1}$) index,
where $i_1$ and $i_{-1}$ denote the cases of $\omega$-index with $\omega=1$ and $\omega=-1$ respectively (cf. Chap. 5 of \cite{lo2002}).
Since we have assume the manifold $M=S^n/ \Gamma$ is odd dimensional in Theorem 1.2,
then $M$ is orientable and we only state the precise index iteration formulae
of orientable closed geodesics in the following. Throughout this section we write $i_1(\gamma)$  as $i(\gamma)$ for short.

For the reader's convenience, we  briefly review some basic materials  in Long's book \cite{lo2002}.

Let $P$ be a symplectic matrix in Sp$(2N-2)$ and $\Omega^{0}(P)$ be the path connected component of its homotopy set $\Omega(P)$ which contains $P$.
Then there is a path $f\in C([0,1],\Omega^{0}(P))$
such that $f(0)=P$ and
 \bea f(1)
&=& N_1(1,1)^{\dm p_-}\,\dm\,I_{2p_0}\,\dm\,N_1(1,-1)^{\dm p_+}\nn\\
& &  \dm\,N_1(-1,1)^{\dm q_-}\,\dm\,(-I_{2q_0})\,\dm\,N_1(-1,-1)^{\dm q_+} \nn\\
&&\dm\,R(\th_1)\,\dm\,\cdots\,\dm\,R(\th_{r'})\,\dm\,R(\th_{r'+1})\,\dm\,\cdots\,\dm\,R(\th_r)\lb{nmf}\\
&&\dm\,N_2(e^{i\aa_{1}},A_{1})\,\dm\,\cdots\,\dm\,N_2(e^{i\aa_{r_{\ast}}},A_{r_{\ast}})\nn\\
& &  \dm\,N_2(e^{i\bb_{1}},B_{1})\,\dm\,\cdots\,\dm\,N_2(e^{i\bb_{r_{0}}},B_{r_{0}})\nn\\
& &\dm\,H(\pm 2)^{\dm h},\nn\eea
where $N_{1}(\lambda,\chi)=\left(
  \begin{array}{ll}
    \lambda\quad \chi\\
   0\quad \lambda\\
  \end{array}
  \right)$ with $\lambda=\pm 1$ and $\chi=0,\ \pm 1$; $H(b)=\left(
  \begin{array}{ll}
    b\quad 0\\
   0\quad b^{-1}\\
  \end{array}
  \right)$ with $b=\pm 2$;
  $R(\theta)= \left(
  \begin{array}{ll}
    \cos\theta\ -\sin\theta\\
   \sin\theta\quad\ \cos\theta\\
  \end{array}
  \right)$ with $\theta\in (0,2\pi)\setminus\{\pi\}$ and we suppose that $\pi<\theta_{j}<2\pi$ iff $1\leq j\leq
r'$;  $$N_{2}(e^{i\aa_{j}},A_{j})=\left(
  \begin{array}{ll}
    R(\aa_{j})\  A_{j}\\
  \ 0\quad\ R(\aa_{j})\\
  \end{array}
  \right)\ \text{and}\ N_{2}(e^{i\bb_{j}},B_{j})=\left(\begin{array}{ll}
    R(\bb_{j})\  B_{j}\\
  \ 0\quad\ R(\bb_{j})\\
  \end{array}
  \right)$$ with $\alpha_{j},\beta_{j}\in (0,2\pi)\setminus\{\pi\}$ are non-trivial and trivial basic normal forms respectively.

Let $\gamma_{0}$ and $\gamma_{1}$ be two symplectic paths in Sp$(2N-2)$ connecting the identity matrix $I$ to $P$ and $f(1)$
satisfying $\gamma_{0}\sim_\omega\gamma_1$. Then it
has been shown that $i_{\omega}(\gamma_{0}^{m})=i_{\omega}(\gamma_{1}^{m})$ for any $\omega\in S^{1}=\{z\in\C\mid|z|=1\}.$
Based on this fact, we always
assume without loss of generality that each $P_{c}$ appearing in the sequel has the form (\ref{nmf}).

 \begin{lemma}\label{orient} {\rm (cf. Theorem 8.3.1 and Chap. 12 of \cite{lo2002})}  Let $c$ be an orientable closed geodesic on
an $N$-dimensional Finsler manifold with its Poincar$\acute{e}$ map $P_c$.  Then, there
exists a continuous symplecitic path  $\Psi$ with $\Psi(0)=I$ and $\Psi(1)=P_c$ such that
\bea i(c^{m})=i(\Psi^m)
&=& m(i(\Psi)+p_-+p_0-r )  -(p_- + p_0+r) - {{1+(-1)^m}\over 2}(q_0+q_+) \nn\\&&
              + 2\sum_{j=1}^r{E}\left(\frac{m\th_j}{2\pi}\right)+ 2\sum_{j=1}^{r_{\ast}}\vf\left(\frac{m\aa_j}{2\pi}\right) - 2r_{\ast},
\lb{preind}\eea
and
\bea
\nu(c^{m})=\nu(\Psi^m)
 &=& \nu(\Psi) + {{1+(-1)^m}\over 2}(q_-+2q_0+q_+) + 2\vs(c,m),    \lb{pnul}\eea
where we denote by
\be \vs(c,m) = \left(r - \sum_{j=1}^r\vf\left(\frac{m\th_j}{2\pi}\right)\right)
             + \left(r_{\ast} - \sum_{j=1}^{r_{\ast}}\vf\left(\frac{m\aa_j}{2\pi}\right)\right)
             + \left(r_0 - \sum_{j=1}^{r_0}\vf\left(\frac{m\bb_j}{2\pi}\right)\right).    \nn\lb{yuxiang}\ee
\end{lemma}

\subsection{A variant of Precise index iteration formulae}
In this subsection, we give a variant of the precise index iteration formulae in section 3.1 which makes them more intuitive and enables us to apply
the Kronecker's approximation theorem to study the multiplicity of non-contractible closed geodesics of the class $[h]$.

To prove Theorem \ref{mainresult}, in the following of this paper,
we always assume that there exists only one non-contractible minimal closed geodesic
$c$ of the class $[h]$ on $S^{2n+1}/ \Gamma$ with a bumpy irreversible metric $F$, which is then just the well known minimal point of the energy functional
$E$ on $\Lambda_{h}M$ satisfying $i(c)=0$. We can suppose that $c=\gamma^t$ for some prime closed geodesic $\gamma$ and $t\in\N$,
then we also have $i(\gamma)=0$ since $\gamma$ is also a local minimal point of the energy functional
$E$.

Now the Morse-type number is given by
\begin{equation*} m_q \equiv M_q(\Lm_hM)=\sum_{m\ge 1}\dim{\ol{C}}_q(E, c^{p(m-1)+1}), \quad
                \forall q\in \N_0. \nn\end{equation*}
Since the Finsler metric F is bumpy, for the Poincar$\acute{e}$ map $P_c$ of $c$, there is a path $f\in C([0,1],\Omega^{0}(P_c))$
such that $f(0)=P_c$ and
 \bea f(1)
&=& R(\th_1)\,\dm\,\cdots\,\dm\,R(\th_k)\,\dm\,N_2(e^{i\aa_{1}},A_{1})\,\dm\,\cdots\,\dm\,N_2(e^{i\aa_{r_{\ast}}},A_{r_{\ast}})\nn\\
& &  \dm\,N_2(e^{i\bb_{1}},B_{1})\,\dm\,\cdots\,\dm\,N_2(e^{i\bb_{r_{0}}},B_{r_{0}})\,\dm\,H(\pm 2)^{\dm h},\nn\eea
where $\frac{\th_{j}}{2\pi}'s$, $\frac{\aa_{j}}{2\pi}'s$ and $\frac{\bb_{j}}{2\pi}'s$ are all in $\Q^c\cap(0,1)$, $k + h + 2r_{\ast}+ 2r_0 = 2n$.
Then by (\ref{preind}) in Lemma \ref{orient} we obtain
\bea i(c^{m})=-mk-k+2\sum_{j=1}^k{E}\left(\frac{m\th_j}{2\pi}\right),
\lb{preind1}\eea
Then we have:
\begin{lemma}
\label{morsebetti} Assuming the existence of only one non-contractible
minimal closed geodesic $c$ of the class $[h]$ on $S^{2n+1}/ \Gamma$ with a bumpy irreversible metric $F$,
where the order of $h$ is $p$ with $p\geq 2$, there hold
\begin{equation} m_{2q+1}=\bar{\beta}_{2q+1}=0\  {\rm and}\ m_{2q}=\bar{\beta}_{2q}, \qquad
           \forall\ q\in \N_0.{\label{beqm}}\end{equation}
\end{lemma}

\Proof We prove in two cases:

{\bf Case 1.} If $i(c^{p+1})-i(c)$ is even.

In this case, by (\ref{preind1}) we have $pk\in2\N$ and $i(c^{p(m-1)+1})-i(c)\in2\Z$ for any $m\in\N$, which implies that
$i(c^{p(m-1)+1})\in2\N_0$ due to $i(c)=0$. Thus by Lemma \ref{Rad1992}(i) and the definition of $m_q$,
note that $i(\gamma)=0$, we obtain $m_{2q+1} = 0$ for
all $q\in\N_0$. Then by (\ref{b.1}) of Lemma \ref{Lm2.3} it yields $m_{2q+1}=\bar{\beta}_{2q+1}= 0$ for all $q\in\N_0$ and
$m_{2q} =\bar{\beta}_{2q}$ follows from the following Morse inequalities:
\bea m_q-m_{q-1}+\cdots+(-1)^qm_0\geq \bar{\beta}_q-\bar{\beta}_{q-1}+\cdots+(-1)^q\bar{\beta}_0, \forall q\in\N_0.\nn\eea

{\bf Case 2.} If $i(c^{p+1})-i(c)$ is odd.

In this case, by (\ref{preind1}) we have $pk\in2\N-1$ and $i(c^{p(m-1)+1})-i(c)\in2\Z$ if and only if $m\in2\N-1$, which implies that
$i(c^{p(m-1)+1})\in2\N_0$ if and only if $m\in2\N-1$ due to $i(c)=0$, then by Lemma \ref{Rad1992}(i) we have
\bea \overline{C}_q( E,c^{p(m-1)+1}) = \left\{\begin{array}{ll}

     \Q, &\quad {\it if}\;m\in 2\N-1\;and\;
                   q=i(c^{p(m-1)+1}),\;  \cr
     0, &\quad {\it otherwise},\\ \end{array}\right. \nn \eea
where we use the fact $i(\gamma)=0$.
Thus by the definition of $m_q$, we obtain $m_{2q+1} = 0$ for
all $q\in\N_0$. Then by (\ref{b.1}) of Lemma \ref{Lm2.3} it yields $m_{2q+1}=\bar{\beta}_{2q+1}= 0$ for all $q\in\N_0$ and
$m_{2q} =\bar{\beta}_{2q}$ follows from the following Morse inequalities:
\bea m_q-m_{q-1}+\cdots+(-1)^qm_0\geq \bar{\beta}_q-\bar{\beta}_{q-1}+\cdots+(-1)^q\bar{\beta}_0, \forall q\in\N_0.\nn\eea
The proof is complete.\hfill$\Box$

Now we prove Theorem \ref{mainresult} for $M=S^{2n+1}/ \Gamma$ with a bumpy reversible Finsler metric $F$:
\begin{theorem}\label{theorem3.1} Let $M=S^{2n+1}/ \Gamma$ and $h$ be a nontrivial element of finite order $p$ in $\pi_1(M)$.
Then every bumpy reversible Finsler metric $F$ on $M$
has at least two distinct non-contractible closed geodesics of the class $[h]$.
\end{theorem}
\Proof Assume that there exists only one non-contractible minimal closed geodesic
$c$ of the class $[h]$ on $M$. When the metric F on $M$ is reversible, the inverse curve $c^{-1}$ of a closed geodesic
$c$ of the class $[h]$ has played the same role in the variational setting of the energy functional $E$ on $\Lambda_h M$ as $c$.
Specially, the $m$-th iterates $c^m$ and $c^{-m}$ have precisely the same Morse indices, nullities and
critical modules. Then by the proof of Lemma \ref{morsebetti}, (\ref{beqm}) also holds for
bumpy reversible Finsler metrics which together with (\ref{b.1}) of Lemma \ref{Lm2.3}
gives $m_0=\bar{\beta}_0=1$. On the other hand, we have
\bea m_0=\sum_{m\ge 1}\dim{\ol{C}}_0(E, c^{\pm(p(m-1)+1)})\geq \dim{\ol{C}}_0(E, c)+\dim{\ol{C}}_0(E, c^{-1})=2,\nn\eea
where we use the fact $i(c^{\pm 1})=i(\gamma)=0$, and Lemma \ref{Rad1992}(i).
This contradiction completes the proof of Theorem \ref{theorem3.1}. \hfill$\Box$

By Theorem \ref{theorem3.1}, we only need to prove Theorem \ref{mainresult} for $M=S^{2n+1}/ \Gamma$ with
a bumpy irreversible Finsler metric $F$ in the following.

\begin{lemma}\label{lemma3.2} Suppose $c$ is the only one non-contractible
minimal closed geodesic $c$ of the class $[h]$ on $S^{2n+1}/ \Gamma$ with a bumpy irreversible metric $F$,
where the order of $h$ is $p$ with $p\geq 2$.
Then there exist an integer $\bar{p}\geq 2$ and
$\hat{\theta}_1$, $\hat{\theta}_2$, \dots, $\hat{\theta}_k$ in $\Q^c$ with $2\leq k\leq2n$ such that $k\bar{p}\in2\N$ and
\bea  \sum_{j=1}^{k}\hat{\theta}_{j}&=&{1\over2}\left(k+{2n\over \bar{p}(n+1)}\right), \label{mean}\\
 i(c^{m})&=& m\left({2n\over \bar{p}(n+1)}\right)+k-2\sum_{j=1}^{k}\left\{{m\hat{\theta}_{j}}\right\},\qquad\forall\ m\ge 1. \lb{iter1}\eea
In addition, $c^{m}$ has contribution to the Morse-type number $\{m_q\mid q\in\N_0\}$ if and only if $m\equiv 1(\mod \bar{p})$.
\end{lemma}
\Proof From (\ref{preind1}), we have \bea \hat{i}(c)=-k+\sum_{j=1}^k\frac{\th_j}{\pi}.
\lb{preind2}\eea

Now we prove in two cases:

{\bf Case 1.} If $i(c^{p+1})-i(c)$ is even.

In this case, by (\ref{preind1}) we have $pk\in2\N$. From (\ref{breident1}) of Theorem \ref{Thm1.1} and Lemma \ref{Rad1992}(i), we have
\be  \frac{(-1)^{i(c)}}{\hat{i}(c)}=
    \frac{p(n+1)}{2n},  \nn\ee
which together with (\ref{preind2}) and the fact $i(c)=0$ yields
\bea  \sum_{j=1}^k\frac{\th_j}{\pi}=k+{2n\over p(n+1)}.\label{mean-1}\eea
Then by (\ref{preind1}) and (\ref{mean-1}), we have \bea i(c^{m})&=&-mk+k+2\sum_{j=1}^k\left[\frac{m\th_j}{2\pi}\right]\nn\\
&=&-mk+k+2\sum_{j=1}^k\left(\frac{m\th_j}{2\pi}-\left\{\frac{m\th_j}{2\pi}\right\}\right)\nn\\
&=&m\left({2n\over p(n+1)}\right)+k-2\sum_{j=1}^k\left\{\frac{m\th_j}{2\pi}\right\}.
\lb{preind-1}\eea
Let $\bar{p}=p$, $\hat{\theta}_j=\frac{\th_j}{2\pi}$ for $j=1,2,\cdots, k$, then $k\bar{p}\in2\N$
and (\ref{mean})-(\ref{iter1}) hold by (\ref{mean-1})-(\ref{preind-1}).
By the proof of Case 1 of Lemma \ref{morsebetti}, we know $c^{m}$ has contribution to the Morse-type
numbers $\{m_q\mid q\in\N_0\}$ if and only if $m\equiv 1(\mod p)$.

{\bf Case 2.} If $i(c^{p+1})-i(c)$ is odd.

In this case, by (\ref{breident1}) of Theorem \ref{Thm1.1} and Lemma \ref{Rad1992}(i), we have
\be  \frac{(-1)^{i(c)}}{\hat{i}(c)}=
    \frac{p(n+1)}{n},  \nn\ee
which together with (\ref{preind2}) and the fact $i(c)=0$ yields
\bea  \sum_{j=1}^k\frac{\th_j}{\pi}=k+{n\over p(n+1)}.\label{mean-2}\eea
Then by (\ref{preind1}) and (\ref{mean-2}), we have \bea i(c^{m})=m\left({n\over p(n+1)}\right)+k-2\sum_{j=1}^k\left\{\frac{m\th_j}{2\pi}\right\}.
\lb{preind-2}\eea
Let $\bar{p}=2p$, $\hat{\theta}_j=\frac{\th_j}{2\pi}$ for $j=1,2,\cdots, k$, then $k\bar{p}\in2\N$
and (\ref{mean})-(\ref{iter1}) hold by (\ref{mean-2})-(\ref{preind-2}).
By the proof of Case 2 of Lemma \ref{morsebetti}, we know $c^{m}$ has contribution to the
Morse-type numbers $\{m_q\mid q\in\N_0\}$ if and only if $m\equiv 1(\mod \bar{p})$.
The proof is complete. \hfill$\Box$

Now we give a variant of the precise index iteration formulae (\ref{iter1}) specially for our purpose. Let  $m=\bar{p}(n+1)l+\bar{p}L+1$
with $l\in\N$ and $L\in\Z$. By (\ref{mean}) and (\ref{iter1}) we obtain
\bea  i(c^{m})
&=& 2nl+k+(\bar{p}L+1){2n\over \bar{p}(n+1) }\nn\\
& & -2\left(\left\{{k\over2}+{(\bar{p}L+1)n\over\bar{p}(n+1)}-\sum_{j=2}^{k}\left\{m\hat{\theta}_{j}\right\}\right\}
    +\sum_{j=2}^{k}\left\{m\hat{\theta_{j}}\right\}\right)  \nn\\
&=& 2nl+2\left[{k\over2}+{(\bar{p}L+1)n\over\bar{p}(n+1)}\right]+2\left\{{k\over2}+{(\bar{p}L+1)n\over\bar{p}(n+1)}\right\}  \nn\\
& & -2\left(\left\{\left\{{k\over2}+{(\bar{p}L+1)n\over\bar{p}(n+1)}\right\}-\sum_{j=2}^{k}\left\{m\hat{\theta}_{j}\right\}\right\}
   +\sum_{j=2}^{k}\left\{m\hat{\theta_{j}}\right\}\right)  \nn\\
&=& 2nl+2\left[Q_{L}\right]+2\left\{Q_{L}\right\}-2\left(\left\{\left\{Q_{L}\right\}-\sum_{j=2}^{k}
  \left\{m\hat{\theta}_{j}\right\}\right\}+\sum_{j=2}^{k}\left\{m\hat{\theta_{j}}\right\}\right),
  \label{evenformulae}\eea
where in the first identity we use the fact $k\bar{p}\in2\N$, in the last identity for notational simplicity, we denote by
\bea Q_{L}={k\over2}+{(\bar{p}L+1)n\over\bar{p}(n+1)}. \label{QL}\eea

Since $\sum_{j=2}^{k}\{m\hat{\theta}_{j}\}\in\Q^{c}$, we obtain by (\ref{evenformulae}) that for $1\leq i\leq k-2$,
\be  i(c^{m})
= \left\{\begin{array}{ll}
2nl+2\left[Q_{L}\right],& \text{iff}\ \sum_{j=2}^{k}\left\{m\hat{\theta}_{j}\right\}\in (0,\left\{Q_{L}\right\}),\\
2nl+2\left[Q_{L}\right]-2i,& \text{iff}\ \sum_{j=2}^{k}\left\{m\hat{\theta}_{j}\right\}\in (i-1+\left\{Q_{L}\right\},i+\left\{Q_{L}\right\}),\\
2nl+2\left[Q_{L}\right]-2(k-1),& \text{iff}\ \sum_{j=2}^{k}\left\{m\hat{\theta}_{j}\right\}\in (k-2+\left\{Q_{L}\right\},k-1).
\end{array}\right.   \lb{evenodddengjia}\ee
Let \bea I_{0}(L)=(0,\left\{Q_{L}\right\}), I_{k-1}(L)=(k-2+\left\{Q_{L}\right\},k-1),\nn\\
 I_{i}(L)=(i-1+\left\{Q_{L}\right\},i+\left\{Q_{L}\right\})\quad \text{for}\quad 1\le i\le k-2.  \label{Qujian}\eea
Then, (\ref{evenodddengjia}) can be stated in short as that for any integers $m=\bar{p}(n+1)l+\bar{p}L+1$ and  $0\leq i\leq k-1,$
\be  i(c^{m}) = 2nl+2\left[Q_{L}\right]-2i\quad \text{if and only if}\quad
    \sum_{j=2}^{k}\left\{m\hat{\theta}_{j}\right\}\in I_{i}(L).   \lb{shortdengjia}\ee

\begin{remark}
\label{transformation} Let $(\tau(1),\tau(2),\dots,\tau(k))$ be an arbitrary permutation of $(1,2,\dots,k)$.
Then, the same conclusion as (\ref{shortdengjia})
with $j$ ranging in $\{\tau(1),\tau(2),\dots,\tau(k-1)\}$ instead is still valid.
\end{remark}

The following lemma will be also needed in the proof of Theorem \ref{mainresult} for $S^{2n+1}/ \Gamma$  in Section 4.
\begin{lemma}\label{bound} Under the assumption of Lemma \ref{lemma3.2}, for any positive integers $l$ and $m$, we have
$$  |i(c^{m})-2nl|>2n\quad \text{holds whenever}\quad  |m-\bar{p}(n+1)l|>2\bar{p}(n+1). $$
\end{lemma}

\Proof From (\ref{iter1}), we have
$$  i(c^{m}) = 2nl+(m-\bar{p}(n+1)l)\cdot{2n\over \bar{p}(n+1)}+k-2\sum_{j=1}^{k}\left\{{m\hat{\theta}_{j}}\right\}, $$
which yields immediately that
\bea  |i(c^{m})-2nl|
&\ge& |m-\bar{p}(n+1)l|\cdot{2n\over \bar{p}(n+1)}-|k-2\sum_{j=1}^{k}\left\{{m\hat{\theta}_{j}}\right\}|  \nn\\
&>& 4n-k \ge 4n-2n = 2n,  \nn\eea
where the fact $k\le 2n$ is used. \hfill$\Box$

\subsection{The system of irrational numbers}
In this subsection, we review some properties of a system of irrational numbers associated
to our proof of Theorem \ref{mainresult}, all the details can be found in section 4 of \cite{LX}.
Let $\aa=\{\alpha_{1}, \alpha_{2}, \ldots, \alpha_{m}\}$ be a set of $m$ irrational numbers. As usual,
we have
\begin{definition}\lb{rank}
The set $\aa$ of irrational numbers is linearly independent over $\Q$, if there do not exist $c_1$,
$c_2$, $\ldots$, $c_m$ in $\Q$ such that $\sum_{j=1}^{m}|c_{j}|> 0$ and
\be  \sum_{j=1}^{m}c_{j}\alpha_{j}\in\Q,  \label{81a}\ee
and is linearly dependent over $\Q$ otherwise. The rank of $\aa$ is defined to be the number of
elements in a maximal linearly independent subset of $\aa$, which we denote by $\rank(\aa)$.
\end{definition}

\begin{lemma}\lb{L3.1}
Let $r=\rank(\aa)$. Then there exist $p_{jl}\in\Z$, $\beta_{l}\in\Q^{c}$ and $\xi_{j}\in \Q$ for $1\leq l\leq r$ and $1\leq j\leq m$ such that
 \begin{equation}
\label{81b}
\alpha_{j}=\sum_{l=1}^{r}p_{jl}\beta_{l}+\xi_{j},\ \forall 1\leq j\leq m.
\end{equation}

\end{lemma}

In order to study the multiplicity of closed geodesics on $S^{2n+1}/ \Gamma$ with a bumpy irreversible metric $F$,
we are particularly interested in the irrational system  $\{\hat{\theta}_{1}$,
$\hat{\theta}_{2}$,\dots,$\hat{\theta}_{k}\}$ with rank $1$ satisfying (\ref{mean}).
Then by Lemma \ref{L3.1}, it can be reduced to the following system
\be   \hat{\theta}_{j}=p_{j}\theta+\xi_{j},\ \forall 1\leq j\leq k,  \label{Qdependent}\ee
with $\theta\in\mathbb{Q}^{c}$, $p_{j}\in \Z\backslash\{0\}$, $\xi_{j}\in\mathbb{Q}\cap[0,1)$
satisfying
\bea
&& p_{1}+p_{2}+\cdots+p_{k}=0, \lb{basic1c}\\
&& \{\xi_{1}+\xi_{2}+\cdots+\xi_{k}\}\in(0,1)\backslash\{{1/2}\}, \lb{basic2c}\eea
where to get $\xi_j\in [0,1)$, if necessary, we can replace $\hat{\th}_j$ and $\xi_j$ by
$\breve{\theta}_{j}=\hat{\theta}_{j}-[\xi_{j}]$ and $\breve{\xi}_{j}=\{\xi_{j}\}$. Now
we verify that (\ref{basic1c})-(\ref{basic2c}) hold. In fact, by (\ref{Qdependent})
we have $\sum_{j=1}^{k}\hat{\theta}_{j}=(\sum_{j=1}^{k}p_{j})\theta+(\sum_{j=1}^{k}\xi_{j})$,
where $\sum_{j=1}^{k}p_{j}\in\Z$, $\theta\in\mathbb{Q}^{c}$, $\sum_{j=1}^{k}\xi_{j}\in\Q$.
On the other hand, by (\ref{mean}) we have $\sum_{j=1}^{k}\hat{\theta}_{j}\in\Q$
and $\{\hat{\theta}_{1}+\hat{\theta}_{2}+\cdots+\hat{\theta}_{k}\}\neq\frac{1}{2}$ or $0$ since ${2n\over \bar{p}(n+1)}\notin \Z$.
Thus (\ref{basic1c})-(\ref{basic2c}) hold.

Take arbitrarily $\eta\in\Q$ and make  the following natural $\eta$-action  to the system (\ref{Qdependent}):
\begin{equation}
\label{eta}
{\eta}(\theta)=\theta+\eta,\ {\eta}(\hat{\theta}_{j})=\hat{\theta}_{j}-\left[\xi_{j}-p_{j}\eta\right]\ \text{and}\ {\eta}({\xi}_{j})
=\left\{\xi_{j}-p_{j}\eta\right\},\ \forall 1\leq j\leq k,
\end{equation}
which is obviously induced by the transformation ${\eta}(\theta)=\theta+\eta$.
Then, we get a new system
\begin{equation}
\label{Qdependent1}
{\eta}(\hat{\theta}_{j})=p_{j}{\eta}(\theta)+{\eta}(\xi_{j}),\ \forall 1\leq j\leq k,
\end{equation}
with
\begin{equation}
\label{invariant1}
\begin{aligned}
\{{\eta}({\xi}_{1})+{\eta}({\xi}_{2})+\cdots+{\eta}({\xi}_{k})\}&=\{\{\xi_{1}-p_{1}\eta\}+\{\xi_{2}-p_{2}\eta\}+\cdots+\{\xi_{k}-p_{k}\eta\}\}\\
&=\{\xi_{1}+\xi_{2}+\cdots+\xi_{k}-(p_{1}+p_{2}+\dots+p_{k})\eta\}\\
&=\{\xi_{1}+\xi_{2}+\cdots+\xi_{k}\},
\end{aligned}
\end{equation}
where the third equality we have used the condition (\ref{basic1c}).
For simplicity of writing, we also denote  the new system (\ref{Qdependent1})  by (\ref{Qdependent})$_{\eta}$ meaning that it comes from
(\ref{Qdependent}) by an $\eta$-action.

For the system (\ref{Qdependent})$_{\eta}$ with $\eta\in\Q$, we divide the set $\{1\leq j\leq k\}$ into the following three parts:
\bea \mathcal{K}_{0}^{+}(\eta)&=&\{1\leq j\leq k\mid \eta(\xi_{j})=0,\ p_{j}>0\}, \nn\\
\mathcal{K}_{0}^{-}(\eta)&=&\{1\leq j\leq k\mid \eta(\xi_{j})=0,\ p_{j}<0\}, \nn\\
\mathcal{K}_{1}(\eta)&=&\{1\leq j\leq k\mid \eta(\xi_{j})\neq0 \}.\lb{K_1}\eea
Denote by $k_{0}^{+}(\eta)$, $k_{0}^{-}(\eta)$ and $k_{1}(\eta)$ the numbers ${}^{\#}\mathcal{K}_{0}^{+}(\eta)$, ${}^{\#}\mathcal{K}_{0}^{-}(\eta)$
and ${}^{\#}\mathcal{K}_{1}(\eta)$ respectively.
For the case of $\eta=0$, we write them for short as $k_{0}^{+}$, $k_{0}^{-}$ and $k_{1}.$
It follows immediately that $$k_{0}^{+}(\eta)+k_{0}^{-}(\eta)+k_{1}(\eta)=k.$$
By (\ref{basic2c}) and (\ref{invariant1}), it is obvious that $k_{1}(\eta)\geq1$ for every $\eta\in\Q$.
\begin{definition}
\label{crucialnumber}
For every $\eta\in\Q$, the {\bf absolute difference number} of (\ref{Qdependent})$_{\eta}$ is
defined to be the non-negative number $|k_{0}^{+}(\eta)-k_{0}^{-}(\eta)|.$
 The {\bf effective difference number} of {(\ref{Qdependent})} is defined by
$$\max\{|k_{0}^{+}(\eta)-k_{0}^{-}(\eta)|\mid \eta\in\Q\}.$$
Two  systems of irrational numbers with rank $1$ are called to be {\bf equivalent}, if their effective difference numbers are the same.
\end{definition}

\begin{remark} By the definition of an $\eta$-action in (\ref{eta}), it can be checked directly that
 $\eta_{1}\circ\eta_{2}=\eta_{1}+\eta_{2}$  for every $\eta_{1}$ and $\eta_{2}$ in $\Q.$ So every system of irrational numbers
 with rank $1$ is equivalent to the one which comes from itself  by an $\eta$-action.
\end{remark}

The following theorem is concerned with the lower estimate on the effective difference number
of (\ref{Qdependent}) and will play a crucial role in our proof of Theorem \ref{mainresult} in Section 4. For the reader's convenience, we give
its proof as an appendix in Section 5.
\begin{theorem}
\label{crucialtheorem}
For every system of irrational numbers (\ref{Qdependent}) satisfying the conditions (\ref{basic1c}) and (\ref{basic2c}),  it holds that
   \begin{equation}
   \label{crucialestimate}
   \max\{|k_{0}^{+}(\eta)-k_{0}^{-}(\eta)|\mid \eta\in\Q\}\geq1.
   \end{equation}
\end{theorem}


\section{Proof of Theorem \ref{mainresult}}

In this section, we prove our main Theorem \ref{mainresult}. By Theorem \ref{theorem3.1},
we only need to prove Theorem \ref{mainresult} for $M=S^{2n+1}/ \Gamma$
with a bumpy irreversible Finsler metric $F$ which is involved in the irrational system
$\{\hat{\theta}_{1}$, $\hat{\theta}_{2}$, \dots, $\hat{\theta}_{k}\}$ with $2\leq k\leq 2n$ satisfying (\ref{mean}).
For sake of readability, we divide it into two cases according to whether $\rank(\hat{\theta}_{1},\hat{\theta}_{2},\dots\hat{\theta}_{k})=1$ or not.
We will give in details the proof for the first case. Based on the well known Kronecker's approximation theorem in Number theory,
the second one can be then proved quite similarly and so we only sketch it.

{\bf Proof of Theorem \ref{mainresult}:} We carry out the proof in two cases.

{\bf Case 1:} $r=\rank(\hat{\theta}_{1},\hat{\theta}_{2},\dots\hat{\theta}_{k})=1.$

As we have mentioned in Section 3.3, the irrational system (\ref{mean}) with $r=1$ can be seen as a special case
of (\ref{Qdependent}) satisfying (\ref{basic1c}) and (\ref{basic2c}).

Since any $\eta$-action with $\eta\in\Q$  to (\ref{Qdependent}), if necessary, does no substantive effect on our following arguments,
by Theorem \ref{crucialtheorem} and Remark \ref{transformation} we can assume without loss of generality that
$$|k_{0}^{+}-k_{0}^{-}|\geq1\ \text{and}\ \mathcal{K}_{1}=\{1,2,\dots, k_{1}\},$$
with $k_{1}\geq1$ due to (\ref{basic2c}), and denote by $\xi_{j}={r_{j}\over q_{j}}$ for $1\leq j\leq k_{1}.$

 Let $\bar{q}=q_{1}q_{2}\cdots q_{k_1}$ and $m_{l}=\bar{p}(n+1)\bar{q}l+1$ with $l\in\N$,
 where $\bar{p}$ is given by Lemma \ref{lemma3.2}. Then by (\ref{Qdependent}) we have
\begin{equation}
\label{sum3}
\begin{aligned}
\sum_{j=2}^{k}\left\{m_{l}\hat{\theta}_{j}\right\}&=\sum_{j=2}^{k_1}\left\{m_{l}\hat{\theta}_{j}\right\}
+\sum_{j=k_{1}+1}^{k}\left\{m_{l}\hat{\theta}_{j}\right\}\\
&=\sum_{j=2}^{k_1}\left\{p_{j}\{m_{l}\theta\}+\xi_{j}\right\}+\sum_{j=k_{1}+1}^{k}\left\{p_{j}\{m_{l}\theta\}\right\},\\
\end{aligned}
\end{equation}
for some $\theta\in\Q^c$. Then the set $\{\{m_{l}\theta\}\mid l\in\N\}$ is dense in $[0,1]$. For every $L\in\Z$,
we introduce the  auxiliary function
\begin{equation}
\label{auxifunction}
f_{L}(x)=\sum_{j=2}^{k_1}\left\{\left\{p_{j}x+\xi_{j}\right\}+\bar{p}L\hat{\theta}_{j}\right\}
+\sum_{j=k_{1}+1}^{k}\left\{\left\{p_{j}x\right\}+\bar{p}L\hat{\theta}_{j}\right\},\ \forall x\in [0,1],
\end{equation}
and denote for simplicity by $f=f_{0}$, which contains only finitely many discontinuous points.

Let $a$ and $b$ in $(0,1)$ be two real numbers  sufficiently close to $0$ and $1$ respectively.
Then,
\begin{equation}
\label{fa}
\begin{aligned}
f(a)=\sum_{j=2}^{k_1}\left\{p_{j}a+\xi_{j}\right\}+\sum_{j=k_{1}+1}^{k}\left\{p_{j}a\right\}
&=\sum_{j=2}^{k_1}\left(p_{j}a+\xi_{j}\right)+\sum_{j\in\mathcal{K}_{0}^{+}}p_{j}a+\sum_{j\in\mathcal{K}_{0}^{-}}\left(1+p_{j}a\right)\\
&=k_{0}^{-}+\sum_{j=2}^{k}p_{j}a+\sum_{j=2}^{k_1}\xi_{j},
\end{aligned}
\end{equation}
and by similar computation,
\begin{equation}
\label{fb}
\begin{aligned}
f(b)&=k_{0}^{+}+\sum_{j=2}^{k}p_{j}(b-1)+\sum_{j=2}^{k_1}\xi_{j}.
\end{aligned}
\end{equation}
It follows by (\ref{fa}) and (\ref{fb}) that
\bea
&&f(a), f(b)\in (0, k-1),\quad [f(a)]\geq 0, \quad [f(b)]\geq 0,\nn\\
|f(b)-f(a)|&=&|k_{0}^{+}-k_{0}^{-}+\sum_{j=2}^{k}p_{j}(b-1-a)|=|k_{0}^{+}-k_{0}^{-}+p_{1}(-b+1+a)|,\label{distance}
\eea
where the second identity we have used $\sum_{j=1}^{k}p_{j}=0.$

\begin{lemma}
 \label{lemma1616a}
Given $\bar{N}\in\N$, for any $a$ and $b$ in $(0,1)$ sufficiently close to $0$ and $1$ respectively, then

{\rm (i)} $f(a)$ and $f(b)$ lie in different intervals of  (\ref{Qujian}) with $L=0$,

{\rm (ii)} $f_{L}(a)$ and $f_{L}(b)$ lie in the same interval of (\ref{Qujian}) for any $1\leq |L|\leq\bar{N}$, including $f_{L}(0).$
\end{lemma}
\Proof (i) By (\ref{distance}) and the assumption, $|f(b)-f(a)|\approx|k_{0}^{+}-k_{0}^{-}|.$
 Here and below, we write $A \approx B$, if $A$ and $B$ can be chosen to be as close to each other as we want.
 Since the length of each interval in (\ref{Qujian}) with $L=0$ is less than or equal to $1$, so $f(a)$ and $f(b)$
 must lie in different ones, provided that  $|k_{0}^{+}-k_{0}^{-}|\geq2$.

 If $|k_{0}^{+}-k_{0}^{-}|=1$,  then $|f(b)-f(a)|\approx1$. For the case of $k=2$,
 since the length of each interval of  (\ref{Qujian}) with $L=0$ is
 less than $1$, (i) follows immediately. The rest case is $k\geq3$, which still contains three subcases.

$1^{\circ}$ If $k_{1}\geq2$, by (\ref{Qdependent})-(\ref{basic1c}), (\ref{mean}) and (\ref{QL}), we have
\bea \left\{\sum_{j=1}^{k_1}\xi_{j}\right\}=\left\{\sum_{j=1}^{k}\xi_{j}\right\}=\left\{\sum_{j=1}^{k}\hat{\theta}_{j}\right\}
=\left\{{k\over2}+{n\over\bar{p}(n+1)}\right\}=\left\{Q_0\right\}.\nn\eea
Then we have
$$\left\{k_{0}^{-}+\sum_{j=2}^{k_1}\xi_{j}\right\}=\left\{-\xi_{1}+\sum_{j=1}^{k_1}\xi_{j}\right\}
=\left\{\left\{\sum_{j=1}^{k_1}\xi_{j}\right\}-\xi_{1}\right\}=\left\{ \{Q_{0}\}-\xi_{1}\right\}.$$
Then we get by (\ref{fa}) that $\{f(a)\}\approx \left\{ \{Q_{0}\}-\xi_{1}\right\}$ or $1$ and $f(a)$ is not equal
to these two numbers since $\sum_{j=2}^{k}p_{j}=-p_1\neq 0$ by (\ref{basic2c}) and the fact that $p_1\neq 0$.
Notice that the dividing points of the intervals in (\ref{Qujian})
with $L=0$ are $$0,\ \{Q_{0}\},\ 1+\{Q_{0}\},\ 2+\{Q_{0}\},\ \dots,\ k-2+\{Q_{0}\},\ k-1.$$ Therefore
$f(a)$ must be an interior point of these intervals.
It then  yields that $f(a)$ and $f(b)$ must lie in two different intervals.

$2^{\circ}$ If $k_{1}=1$ and $k_{0}^{-}\geq1$, then  $f(a)\approx k_{0}^{-}$ is also an interior point and (i) follows.

$3^{\circ}$ If $k_{1}=1$ and $k_{0}^{-}=0$, then $f(a)=\sum_{j=2}^{k}p_{j}a$ lies in the first interval
whose length is $\{Q_{0}\}<1$ and so $f(b)$ must lie in another one.

(ii) It can be checked directly that $\lim_{a\to0}f_{L}(a)=\lim_{b\to1}f_{L}(b)=f_{L}(0)\in\mathbb{Q}^{c},$
since $\xi_{j}\in\Q$ for $1\leq j\leq k$ and $\sum_{j=2}^{k}\bar{p}L\hat{\theta}_{j}\in\mathbb{Q}^{c}$ by (\ref{mean}).
But the dividing points of these intervals in (\ref{Qujian}) with $1\leq|L|\leq\bar{N}$ are finitely many rational numbers,
so $f_{L}(0)$ is an interior point of these intervals and (ii) follows. \hfill$\Box$

Notice that $f$ contains only finitely many discontinuous points on $(0,1)$.
Without loss of generality, we assume $a$ and $b$ to be two continuous points of $f$ and choose $l_{1}$, $l_{2}\in\N$ with $l_{2}-l_{1}$
sufficiently large such that $\{m_{l_1}\theta\}\approx a$ and $\{m_{l_2}\theta\}\approx b.$
Then by (\ref{sum3}), (\ref{auxifunction}) and (i) of Lemma \ref{lemma1616a}, we get
$\sum_{j=2}^{k}\left\{m_{l_1}\hat{\theta}_{j}\right\}$ and $\sum_{j=2}^{k}\left\{m_{l_2}\hat{\theta}_{j}\right\}$
lie in different intervals of (\ref{Qujian}) with $L=0$.
Suppose that $$\sum_{j=2}^{k}\left\{m_{l_1}\hat{\theta}_{j}\right\}\in I_{i'}\ \text{and}\ \sum_{j=2}^{k}\left\{m_{l_2}\hat{\theta}_{j}\right\}\in I_{i''},$$
with $\{i',i''\}\subseteq\{0,1,2,\dots,k-1\}$ and $i'\neq i''$. By (\ref{shortdengjia}) we have
$i(c^{m_{l_{1}}})=2n\bar{q}l_{1}+2[Q_{0}]-2i'$ and
\begin{equation}
\label{73a}
i(c^{m_{l_2}})=2n\bar{q}l_{2}+2[Q_{0}]-2i''.
\end{equation}
Since $2n\mid(2n\bar{q}l_{1}+2[Q_{0}]-2i'')$ if and only if $2n\mid(2n\bar{q}l_{2}+2[Q_{0}]-2i'')$, we get by (\ref{b.1})
with $n$ there replaced by $2n+1$ that
\bea \bar{\beta}_{2n\bar{q}l_{1}+2[Q_{0}]-2i''}=\bar{\beta}_{2n\bar{q}l_{2}+2[Q_{0}]-2i''}\equiv\beta.\lb{BET}\eea
 Take $\bar{N}>2(n+1)$ in (ii) of Lemma \ref{lemma1616a} and observe that\bea
 |2[Q_{0}]-2i''|&=&\left|2\left[{k\over2}+{n\over\bar{p}(n+1)}\right]-2i''\right|\nn\\&\leq&
 \max {\left\{2\left[{k\over2}+{n\over\bar{p}(n+1)}\right], 2(k-1)-2\left[{k\over2}+{n\over\bar{p}(n+1)}\right]\right\}}\nn
 \\&\leq& k\leq2n.\lb{Q0}\eea By Lemma \ref{morsebetti} and Lemma \ref{bound},
 there exist $L_{i}\in\Z$ with $1\leq |L_{i}|\leq \bar{N}$ and  $1\leq i\leq\beta$ such that
$$i(c^{m_{l_{1}}+\bar{p}L_{i}})=2n\bar{q}l_{1}+2[Q_{0}]-2i''.$$
In fact, by (\ref{b.1}), $\bar{\beta}_q\neq0$ whenever $q$ is even, then by Lemma \ref{morsebetti}, for the even integer
$2n\bar{q}l_{1}+2[Q_{0}]-2i''$, $\exists$ an integer $\bar{m}$ such that $i(c^{\bar{m}})=2n\bar{q}l_{1}+2[Q_{0}]-2i''$.
From Lemma 3.3, we have $\bar{m}\equiv 1(\mod \bar{p})$. By definition, $m_{l_1}=\bar{p}(n+1)\bar{q}l_1+1\equiv 1(\mod \bar{p})$,
then $\exists$ $L_i\in\Z$ such that $\bar{m}=m_{l_1}+\bar{p}L_i$. But by (\ref{Q0}), we obtain $0\leq |L_{i}|\leq \bar{N}$ by Lemma
\ref{bound}. From $i'\neq i''$, we obtain $1\leq |L_{i}|\leq \bar{N}$.

By (\ref{sum3})-(\ref{auxifunction}), we have $f_{L_i}(a)\approx\sum_{j=2}^{k}\left\{(m_{l_{1}}+\bar{p}L_{i})\hat{\theta}_{j}\right\}$
and $f_{L_i}(b)\approx\sum_{j=2}^{k}\left\{(m_{l_{2}}+\bar{p}L_{i})\hat{\theta}_{j}\right\}$
since $\{m_{l_1}\theta\}\approx a$ and $\{m_{l_2}\theta\}\approx b$.
Thus $\sum_{j=2}^{k}\left\{(m_{l_{1}}+\bar{p}L_{i})\hat{\theta}_{j}\right\}$
and $\sum_{j=2}^{k}\left\{(m_{l_{2}}+\bar{p}L_{i})\hat{\theta}_{j}\right\}$ are in the same interval of (\ref{Qujian})
with $1\leq |L_{i}|\leq\bar{N}$ by (ii) of Lemma \ref{lemma1616a}, we get by (\ref{shortdengjia}) that
\begin{equation}
\label{73b}
i(c^{m_{l_{2}}+\bar{p}L_{i}})=2n\bar{q}l_{2}+2[Q_{0}]-2i'',\ \forall 1\leq i\leq \beta
\end{equation}
By (\ref{BET}), it yields $\beta\equiv\bar{\beta}_{2n\bar{q}l_{2}+2[Q_{0}]-2i''}$. Combining (\ref{73a}) with (\ref{73b}), there
are at least $\beta+1$ iterates of $c$ possessing Morse index $\hat{k}\equiv2n\bar{q}l_{2}+2[Q_{0}]-2i''$.
By the bumpy condition, they all contribute to the Morse type number $m_{\hat{k}}$.
This proves $\beta=\bar{\beta}_{\hat{k}}=m_{\hat{k}}\geq \beta+1$ which is obviously absurd.

{\bf Case 2:} $r=\rank(\hat{\theta}_{1},\hat{\theta}_{2},\dots\hat{\theta}_{k})\geq2.$

By Lemma \ref{L3.1},  there are $p_{jl}\in\Z$, $\theta_{k_l}\in\Q^{c}$ and $\xi_{j}\in\Q$ with $1\leq l\leq r$ and $1\leq j\leq k$ such that
\begin{equation}
\label{system82a}
\hat{\theta}_{j}=\sum_{l=1}^{r}p_{jl}\theta_{k_l}+\xi_{j},\ \forall 1\leq j\leq k.
\end{equation}
Moreover,  $\theta_{k_1}$, $\theta_{k_2}$ \dots, $\theta_{k_r}$ are linearly independent over $\Q$. Due to
(\ref{mean}), it follows
\begin{equation}
\label{20150923a}
\sum_{j=1}^{k}p_{jl}=0,\ \forall 1\leq l\leq r.
\end{equation}
Then we have \bea
\{\xi_{1}+\xi_{2}+\cdots+\xi_{k}\}\in(0,1)\backslash\{{1/2}\}. \lb{4.5*}\eea
In fact, by (\ref{system82a}) and (\ref{20150923a})
we have $\sum_{j=1}^{k}\hat{\theta}_{j}=\sum_{l=1}^{r}(\sum_{j=1}^{k}p_{jl})\theta_{k_l}+\sum_{j=1}^{k}\xi_{j}=\sum_{j=1}^{k}\xi_{j}$,
and from (\ref{mean}) we have $\{\hat{\theta}_{1}+\hat{\theta}_{2}+\cdots+\hat{\theta}_{k}\}\neq\frac{1}{2}$ or $0$
since ${2n\over \bar{p}(n+1)}\notin \Z$.
Thus (\ref{4.5*}) holds.

Our basic idea for proving Case 2 is to construct an irrational system  with rank $1$  associated to (\ref{system82a}),
which plays the essential role in our sequel arguments  due to  the following result.\\
{\bf Kronecker's approximation theorem} (cf. Theorem 7.10 in \cite{Apostol}): {\it If $\theta_{1}$,
$\theta_{2}$, \dots, $\theta_{r}$ are linearly independent over $\Q$,
then the set of all vectors of the form $(\{m\theta_{1}\},\{m\theta_{2}\},\dots,\{m\theta_{r}\})$ for all $m\in\N$ is dense in
$$[0,1]^{r}=\underbrace{[0,1]\times[0,1]\times\cdots\times[0,1]}_{r}.$$}
\begin{lemma}
\label{lemma82a}
For the integers $p_{jl}$'s in (\ref{system82a})-(\ref{20150923a}), there are $s_{2},\ s_{3}, \dots,\ s_{r}\in\Z$ such that
\begin{equation}
\label{82a}
p_{j1}+\sum_{l=2}^{r}s_{l}p_{jl}\in\Z\backslash\{0\},\ \forall 1\leq j\leq k,
\end{equation}
\end{lemma}
\Proof  Let $J_{0}=\{1\leq j\leq k\mid p_{j1}=0\}.$ If $J_{0}=\emptyset$, we need only take $s_{2}=s_{3}=\cdots=s_{r}=0.$
 If $J_{0}\neq\emptyset$,  we claim that $(p_{j2},p_{j3},\dots,p_{jr})\neq(0,0,\dots,0)$ for each $j\in J_{0}$.
 Otherwise,  then (\ref{system82a}) yields that
$\hat{\theta}_{j}=\xi_{j}\in\Q,$
which contradicts to $\hat{\theta}_{j}\in\Q^{c}$. So the set
\begin{equation*}
\label{82b}
X_{j}\equiv\left\{(x_{2},x_{3},\dots,x_{r})\mid p_{j2}x_{2}+p_{j3}x_{3}+\dots+p_{jr}x_{r}=0\right\},
\end{equation*}
 is a subspace of dimension $r-2$ in $\R^{r-1}$ which yields that $X=\cup_{j\in J_{0}}X_{j}$ is a proper subset of $\R^{r-1}$.
 Pick up an arbitrary integral point $(\bar{s}_{2},\bar{s}_{3},\dots,\bar{s}_{r})\in\R^{r-1}\backslash X$.
Then for every $\bar{N}\in\N$ we have
\begin{equation}
\label{20150922a}
|p_{j1}+\sum_{l=2}^{r}\bar{N}\bar{s}_{l}p_{jl}|
=\left\{\begin{array}{ll}
\bar{N}|\sum_{l=2}^{r}\bar{s}_{l}p_{jl}|\neq0,& \text{if}\ j\in J_{0},\\
|p_{j1}|\neq0,& \text{if}\ j\notin J_{0}\ \text{and}\ \sum_{l=2}^{r}\bar{s}_{l}p_{jl}=0,\\
|p_{j1}+\bar{N}\sum_{l=2}^{r}\bar{s}_{l}p_{jl}|,& \text{if}\ j\notin J_{0}\ \text{and}\ \sum_{l=2}^{r}\bar{s}_{l}p_{jl}\neq0.
\end{array}
\right.
\end{equation}
For the third case in the righthand side of (\ref{20150922a}), we can take $\bar{N}\in\N$ sufficiently large so that
$|p_{j1}+\bar{N}\sum_{l=2}^{r}\bar{s}_{l}p_{jl}|\neq0$ for all these $j$'s therein. Finally let $s_{l}=\bar{N}\bar{s}_{l}$
and (\ref{82a}) follows. \hfill$\Box$

Let $\td{p}_{j1}=p_{j1}+\sum_{l=2}^{r}s_{l}p_{jl}\in\Z\backslash\{0\}$, $\td{p}_{jl}=p_{jl}$ if $2\leq l\leq r$.
By Lemma \ref{lemma82a}, we can make the change of variables
$\td{\theta}_{k_1}=\theta_{k_1}\ \text{and}\ \td{\theta}_{k_l}=\theta_{k_l}-s_{l}\theta_{k_1}\ \text{for}\ 2\leq l\leq r.$
Then the system (\ref{system82a}) is transformed to
\begin{equation}
\label{system82b}
\hat{\theta}_{j}=\sum_{l=1}^{r}\td{p}_{jl}\td{\theta}_{k_l}+\xi_{j},\ \forall 1\leq j\leq k,
\end{equation}
and by (\ref{20150923a}) we have
\begin{equation*}
\label{basic2ccc}
\sum_{j=1}^{k}\td{p}_{j1}=\sum_{j=1}^{k}p_{j1}+\sum_{j=1}^{k}\sum_{l=2}^{r}s_{l}p_{jl}=0+\sum_{l=2}^{r}s_{l}\left(\sum_{j=1}^{k}p_{jl}\right)=0.
\end{equation*}
Since  $\theta_{k_1}$, $\theta_{k_2}$, \dots, $\theta_{k_r}$ are linearly independent over $\Q$,
so do  $\td{\theta}_{k_1}$, $\td{\theta}_{k_2}$, \dots, $\td{\theta}_{k_r}.$

Consider the following irrational system  with rank $1$ associated to (\ref{system82b})
\begin{equation}
\label{system82c}
\hat{\alpha}_{j}=\td{p}_{j1}\td{\theta}_{k_1}+\xi_{j},\ \forall 1\leq j\leq k.
\end{equation}
By Theorem \ref{crucialtheorem} and the properties of $\xi_{j}$ in (\ref{system82a}) and (\ref{4.5*}),
without loss of generality we can assume for (\ref{system82c}) that $|\td{k}_{0}^{+}-\td{k}_{0}^{-}|\geq1$ and
denote the corresponding integer set in (\ref{K_1}) by
$\td{\mathcal{K}}_{1}(0)=\{1,2,\dots, \td{k}_{1}\}$, and denote $\xi_{j}$ by $\xi_{j}={r_{j}\over q_{j}}$
with $(r_j, q_j)=1$ for $1\leq j\leq \td{k}_{1}.$

 Let $\td{q}=q_{1}q_{2}\cdots q_{\td{k}_{1}}$ and $\td{m}_{l}=\bar{p}(n+1)\td{q}l+1$ for $l\in\N$, where $\bar{p}$ is given by Lemma \ref{lemma3.2}.
 Then, we get by (\ref{system82b}) that
\begin{equation}
\label{sum33}
\begin{aligned}
\sum_{j=2}^{k}\left\{\td{m}_{l}\hat{\theta}_{j}\right\}
&=\sum_{j=2}^{\td{k}_1}\left\{\td{m}_{l}\hat{\theta}_{j}\right\}+\sum_{j=\td{k}_{1}+1}^{k}\left\{\td{m}_{l}\hat{\theta}_{j}\right\}\\
&=\sum_{j=2}^{\td{k}_1}\left\{\sum_{l=1}^{r}\td{p}_{jl}\{\td{m}_{l}\td{\theta}_{k_l}\}+\xi_{j}\right\}
+\sum_{j=\td{k}_{1}+1}^{k}\left\{\sum_{l=1}^{r}\td{p}_{jl}\{\td{m}_{l}\td{\theta}_{k_l}\}\right\},\\
\end{aligned}
\end{equation}
where note that $\xi_{j}=0$ when $\td{k}_1+1\leq j\leq k$.
By Kronecker's approximation theorem,  the set
$\{(\{\td{m}_{l}\td{\theta}_{k_1}\},\{\td{m}_{l}\td{\theta}_{k_2}\},\dots,\{\td{m}_{l}\td{\theta}_{k_r}\})\mid l\in\N\}$ is dense in $[0,1]^{r}$.
For every $L\in\Z$, similarly to (\ref{auxifunction}) we can introduce the  auxiliary multi-variable function on $[0,1]^{r}$,
$$g_{L}(x_{1},x_{2},\dots,x_{r})=\sum_{j=2}^{\td{k}_1}\left\{\sum_{l=1}^{r}\td{p}_{jl}x_{l}+\xi_{j}+\bar{p}L\hat{\theta}_{j}\right\}
+\sum_{j=\td{k}_{1}+1}^{k}\left\{\sum_{l=1}^{r}\td{p}_{jl}x_{l}+\bar{p}L\hat{\theta}_{j}\right\},$$
and denote for simplicity by $g=g_{0}.$ Similarly as before, we have
 \begin{lemma}
\label{lemma1616b}
Given $\bar{N}\in\N$, let $(a_{1},a_{2},\dots, a_{r})$ and $(b_{1},b_{2},\dots, b_{r})$ in $(0,1)^{r}$
be sufficiently close to $(0,0,0,\dots,0)$ and $(1,0,0,\dots,0)$ respectively,
then

{\rm (i)} $g(a_{1},a_{2},\dots, a_{r})$ and $g(b_{1},b_{2},\dots, b_{r})$ lie in different intervals of  (\ref{Qujian}) with $L=0$,
if we further require $\frac{a_2+\cdots+a_r}{a_1}$ and $\frac{b_2+\cdots+b_r}{1-b_1}$ are sufficiently small.

 {\rm (ii)} $g_{L}(a_{1},a_{2},\dots, a_{r})$ and $g_{L}(b_{1},b_{2},\dots, b_{r})$
lie in the same interval of (\ref{Qujian}) for any $1\leq |L|\leq\bar{N}$, including $g_{L}(0,0,\dots, 0).$
\end{lemma}
\Proof
 (i) Since $a_1$, $a_2$, \dots, $a_r$ (resp. $b_1$, $b_2$, \dots, $b_r$) are independent,
 we can select them by such a way that the decimal functions in
 $g(a_{1},a_{2},\dots, a_{r})$ and $g(b_{1},b_{2},\dots, b_{r})$ are mainly determined by $a_1$ and $b_1$ respectively.
 For instance, this can be realized by requiring $a_l$ (resp. $b_l$) with $2\leq l\leq r$ to be much smaller than $a_1$ (resp. $1-b_1$).
 The rest proof is then similar as that in Lemma \ref{lemma1616a}-(i), with $g$ in stead of $f$ therein.

 (ii) follows the same line as  Lemma \ref{lemma1616a}-(ii)
 and do not need such choices on $\frac{a_2+\cdots+a_r}{a_1}$ and $\frac{b_2+\cdots+b_r}{1-b_1}$ as above.

Due to Lemma \ref{lemma1616b}, the rest proof is then almost word by word as that in Case 1 and so we omit the tedious details.
We complete the proof of Theorem 1.2.\hfill$\Box$

\setcounter{equation}{0}
\section{Appendix }
For the reader's convenience, we give the proof of Theorem \ref{crucialtheorem} as an appendix in this section.

\begin{lemma}
\label{equivalentsystem}
Assume that
\begin{equation}
\label{20150729a}
\left\{\begin{array}{ll}
\hat{\theta}_{j}=p_{j}\theta+\xi_{j},& \forall 1\leq j\leq k-1,\\
\hat{\theta}_{k}=p_{k}\theta,\\
\end{array}
\right.
\end{equation}
with $\sum_{j=1}^{k}p_{k}=0$ and $\left\{\sum_{j=1}^{k-1}\xi_{k}\right\}\in(0,1)\backslash\{1/2\}.$

Then, (\ref{20150729a}) is equivalent to
\begin{equation}
\label{20150729b}
\left\{\begin{array}{ll}
\hat{\theta}_{j}=p_{j}\theta+\xi_{j},& \forall\ 1\leq j\leq k-1,\\
\hat{\theta}_{k,l}=\text{sgn}(p_{k})\theta+{l\over |p_{k}|},& \forall\ 0\leq l\leq |p_{k}|-1,\\
\end{array}
\right.
\end{equation}
 where as usual we define $\sgn(a)=\pm 1$ for $a\in\R\bs\{0\}$ when $\pm a>0$.
\end{lemma}
\Proof Take $\eta\in\Q$ arbitrarily and recall the definition of $\eta$-action in (\ref{eta}).
Then the equation $\hat{\theta}_{k}=p_{k}\theta$ contributes $\text{sgn}(p_{k})$
to the absolute difference number of (\ref{20150729a})$_{\eta}$ if and only if
$$\eta(0)=\{0-p_{k}\eta\}=\{-p_{k}\eta\}=0,$$ that is $\eta\in\Z_{|p_{k}|}$,
which is also the sufficient and necessary condition such that the equations
 $$\hat{\theta}_{k,l}=\text{sgn}(p_{k})\theta+{l\over |p_{k}|},\ \forall\ 0\leq l\leq |p_{k}|-1,$$
 contribute $\text{sgn}(p_{k})$ to the absolute difference number of (\ref{20150729b})$_{\eta}$.
 Since the other equations with $1\leq j\leq k-1$  in (\ref{20150729a})
 and (\ref{20150729b}) are the same, so do their contributions to the absolute difference numbers
 of (\ref{20150729a})$_{\eta}$ and (\ref{20150729b})$_{\eta}$.
 As a result, the absolute difference numbers of (\ref{20150729a})$_{\eta}$ and (\ref{20150729b})$_{\eta}$ are equal for any $\eta\in\Q$
 which yields that the effective difference numbers of (\ref{20150729a}) and (\ref{20150729b})
 are the same and so they are equivalent. \hfill$\Box$
\begin{remark}
\label{beautifulchange}
For the system (\ref{20150729b}), we have
$$\left\{\sum_{j=1}^{k-1}\xi_{j}+\sum_{l=0}^{|p_{k}|-1}{l\over |p_{k}|}\right\}
=\left\{\begin{array}{ll}
\left\{\sum_{j=1}^{k-1}\xi_{j}\right\},& if\ p_{k}\ is\ odd,\\
\left\{\sum_{j=1}^{k-1}\xi_{j}+{1\over2}\right\},& if\ p_{k}\ is\ even.
\end{array}
\right.$$
By the assumption of $\left\{\sum_{j=1}^{k-1}\xi_{j}\right\}\in(0,1)\backslash\{1/2\}$, it follows that
$$\left\{\sum_{j=1}^{k-1}\xi_{j}+\sum_{l=0}^{|p_{k}|-1}{l\over |p_{k}|}\right\}\in(0,1)\backslash\{1/2\}.$$
\end{remark}

\begin{lemma}
\label{cutoff}
If there exist $1\leq j'<j''\leq k$ satisfying that $p_{j'}\cdot p_{j''}=-1$ and
$\left\{{\xi_{j'}}+ {\xi_{j''}}\right\}=0$ in
\begin{equation}
\label{Qdependent1029}
\hat{\theta}_{j}=p_{j}\theta+\xi_{j},\ \forall 1\leq j\leq k,
 \end{equation}
 then (\ref{Qdependent1029}) is equivalent to the system
\begin{equation}
\label{cutoffequations}
\hat{\theta}_{j}=p_{j}\theta+\xi_{j},\ \forall j\in\{1,2,\dots,k\}\backslash\{j',j''\}.
\end{equation}
\end{lemma}
\Proof Assume without loss of generality that $p_{j'}=-p_{j''}=1$ and take $\eta\in\Q$
arbitrarily. Then by (\ref{eta}) and the given condition, we have
$$
\left\{{\eta}(\xi_{j'})+{\eta}(\xi_{j''})\right\}=\left\{ \{\xi_{j'}-\eta\} + \{\xi_{j''}+\eta\} \right\}=\{\xi_{j'}+\xi_{j''}\}=0.
$$
Thus, ${\eta}({\xi_{j'}})=0$ if and only if ${\eta}({\xi_{j''}})=0,$ that is, $j'\in\mathcal{K}_{0}^{+}({\eta})$ if and only if
$j''\in\mathcal{K}_{0}^{-}({\eta}).$ As a result, $p_{j'}$ and $p_{j''}$ together contribute nothing to the absolute difference number
of (\ref{Qdependent1029})$_\eta$ for any $\eta\in\Q$. It then follows immediately that (\ref{Qdependent1029}) is equivalent to
(\ref{cutoffequations}). \hfill$\Box$

{\bf Proof of Theorem \ref{crucialtheorem}:} We carry out the proof with two steps.

{\bf Step 1:} First, letting $\eta_{k}={\xi_{k}\over p_{k}}$ and making ${\eta_{k}}$-action to the original system (\ref{Qdependent}),
we obtain by (\ref{eta}) that
\begin{equation}
\label{729change1}
\left\{\begin{array}{ll}
{\eta_{k}}(\hat{\theta}_{j})=p_{j}{\eta_{k}}(\theta)+{\eta_{k}}(\xi_{j}),\ \forall 1\leq j\leq k-1,\\
{\eta_{k}}(\hat{\theta}_{k})=p_{k}{\eta_{k}}(\theta).
\end{array}
\right.
\end{equation}
Then by Lemma \ref{equivalentsystem}, the system (\ref{729change1}) is equivalent to
\begin{equation}
\label{729change11}
\left\{\begin{array}{ll}
{\eta_{k}}(\hat{\theta}_{j})=p_{j}{\eta_{k}}(\theta)+{\eta_{k}}(\xi_{j}),& \forall\ 1\leq j\leq k-1,\\
\hat{\theta}_{k,l'}=\text{sgn}(p_{k}){\eta_{k}}(\theta)+{l'\over |p_{k}|},& \forall\ 0\leq l'\leq |p_{k}|-1,\\
\end{array}
\right.
\end{equation}

Secondly, taking $\eta_{k-1}\in\Q$ such that ${\eta_{k-1}}\circ{\eta_{k}}(\xi_{k-1})=0$
and making ${\eta_{k-1}}$-action to the system (\ref{729change11}), we get
\begin{equation}
\label{729change2}
\left\{\begin{array}{ll}
{\eta_{k-1}}\circ{\eta_{k}}(\hat{\theta}_{j})=p_{j}{\eta_{k-1}}\circ{\eta_{k}}(\theta)+{\eta_{k-1}}\circ{\eta_{k}}(\xi_{j}),& \forall\ 1\leq j\leq k-2,\\
{\eta_{k-1}}\circ{\eta_{k}}(\hat{\theta}_{k-1})=p_{k-1}{\eta_{k-1}}\circ{\eta_{k}}(\theta),& \\
{\eta_{k-1}}\circ(\hat{\theta}_{k,l'})=\text{sgn}(p_{k}){\eta_{k-1}}\circ{\eta_{k}}(\theta)+{\eta_{k-1}}({l'\over |p_{k}|}),& \forall\ 0\leq l'\leq |p_{k}|-1.\\
\end{array}
\right.
\end{equation}
 Again by Lemma \ref{equivalentsystem}, the system (\ref{729change2}) is equivalent to
\begin{equation}
\label{729change22}
\left\{\begin{array}{ll}
{\eta_{k-1}}\circ{\eta_{k}}(\hat{\theta}_{j})=p_{j}{\eta_{k-1}}\circ{\eta_{k}}(\theta)+\eta_{k-1}\circ{\eta_{k}}(\xi_{j}),& \forall\ 1\leq j\leq k-2,\\
\hat{\theta}_{k-1,l''}=\text{sgn}(p_{k-1}){\eta_{k-1}}\circ{\eta_{k}}(\theta)+{l''\over |p_{k-1}|},& \forall\ 0\leq l''\leq |p_{k-1}|-1,\\
{\eta_{k-1}}(\hat{\theta}_{k,l'})=\text{sgn}(p_{k}){\eta_{k-1}}\circ{\eta_{k}}(\theta)+{\eta_{k-1}}({l'\over |p_{k}|}),& \forall\ 0\leq l'\leq |p_{k}|-1,\\
\end{array}
\right.
\end{equation}

Repeating the above procedure for the rest equations with $j={k-2}$, ${k-3}$, $\cdots$, $2$, $1$ one at a time in order,
we can finally get a system equivalent to the original system (\ref{Qdependent}) which can be written in a simple form such as
\begin{equation}
\label{finalsystem}
\hat{\alpha}_{jl}=\text{sgn}(p_{j})\alpha+\xi_{jl},\ \forall\  1\leq j\leq k\ \text{and}\ 0\leq l\leq |p_{j}|-1,
\end{equation}
with $\alpha\in\Q^{c}$ and $\xi_{jl}\in\Q\cap[0,1)$. Moreover,
by (\ref{invariant1}) and Remark \ref{beautifulchange} we have
\begin{equation}
\label{incredible}
\left\{\sum_{j=1}^{k}\sum_{l=0}^{|p_{j}|-1}\xi_{jl}\right\}\in(0,1)\backslash\{{1/2}\}.
\end{equation}

{\bf Step 2:} We can cut off all the superfluous equations of the system (\ref{finalsystem}), if there are such pairs as that in Lemma \ref{cutoff}.
That is, (\ref{finalsystem}) is equivalent to some a system
\begin{equation}
\label{finalsystem1}
\hat{\theta}_{i}^{\prime}=p_{i}^{\prime}\alpha+\xi_{i}^{\prime},\ \forall\  1\leq i\leq \bar{k},
\end{equation}
with $|p_{i}^{\prime}|=1$, $\sum_{i=1}^{\bar{k}}p_{i}^{\prime}=0$ and
\begin{equation}
\label{basic2cc}
\left\{\sum_{i=1}^{\bar{k}}\xi_{i}^{\prime}\right\}\in(0,1)\backslash\{1/2\}.
\end{equation}
Here notice that $\bar{k}\geq1$ is ensured by the condition (\ref{basic2cc}).

Since all the superfluous equations are cut off, it follows that ${\bar{k}}_{0}^{+}\cdot {\bar{k}}_{0}^{-}=0.$
Assume without loss of generality that
${\bar{k}}_{0}^{+}={\bar{k}}_{0}^{-}=0,$  otherwise we have nothing to do. Since $\sum_{i=1}^{\bar{k}}p_{i}^{\prime}=0$,
we get $${}^{\#}\{1\leq i\leq \bar{k}\mid p_{i}^{\prime}=1\}={}^{\#}\{1\leq i\leq \bar{k}\mid p_{i}^{\prime}=-1\}.$$
Take arbitrarily out $i_{1}\in \{1\leq i\leq \bar{k}\mid p_{i}^{\prime}=1\}.$ Let $\bar{\eta}=\xi_{i_{1}}'$
and make the $\bar{\eta}$-action  to (\ref{finalsystem1}).
Then it follows immediately that  $\bar{k}_{0}^{+}(\bar{\eta})\geq1$.
Recalling again that all the superfluous equations have been cut off at the beginning of Step 2,
we obtain $\bar{\eta}(\xi_{i}')=\{\xi_{i_{1}}'+\xi_{i}'\}\neq0$ for
every $i\in\{1\leq i\leq \bar{k}\mid p_{i}^{\prime}=-1\}$ which yields $\bar{k}_{0}^{-}(\bar{\eta})=0.$  As a result, we get
$$\max\{|\bar{k}_{0}^{+}({\eta})-\bar{k}_{0}^{-}({\eta})\mid\eta\in\Q\}\geq|\bar{k}_{0}^{+}(\bar{\eta})
-\bar{k}_{0}^{-}(\bar{\eta})|=\bar{k}_{0}^{+}(\bar{\eta})\geq1.$$
Since the original system (\ref{Qdependent}) is equivalent to (\ref{finalsystem1}), the estimate (\ref{crucialestimate}) follows immediately. \hfill$\Box$

\end{document}